\documentclass[12pt,reqno]{amsart}
\usepackage{amsmath, amsthm, amssymb, amscd}
\usepackage{graphicx, texdraw, lscape}

\input texdraw
\input xy
\xyoption{all}

\newcommand{\Conv}{\mathrm{Conv}}
\newcommand{\Int}{\mathrm{Int}}
\newcommand{\medno}{\medskip\noindent}
\newcommand{\rank}{\mathrm{rank}}

\swapnumbers
\theoremstyle{plain}  
\newtheorem{theorem}[subsection]{Theorem}
\newtheorem{proposition}[subsection]{Proposition}
\newtheorem{lemma}[subsection]{Lemma}
\newtheorem{corollary}[subsection]{Corollary}

\begin{document}
\topmargin-0.75truecm

\title[Noncrystallographic dominant regions]%
{Dominant regions in noncrystallographic \\
hyperplane arrangements }

\author{Yu Chen}
\address{Department of Mathematics \\
Idaho State University \\
Pocatello, ID 83209-8085} \email{chenyu@isu.edu}
\author{Cathy Kriloff}
\address{Department of Mathematics \\
Idaho State University \\
Pocatello, ID 83209-8085} \email{krilcath@isu.edu}

\thanks{
Research of the second author supported in part by
National Security Agency grant MDA904-03-1-0093.}

\maketitle

\begin{abstract}
For a crystallographic root system, dominant regions in the
Catalan hyperplane arrangement are in bijection with antichains in
a partial order on the positive roots.  For a noncrystallographic
root system, the analogous arrangement and regions have
importance in the representation theory of an associated graded
Hecke algebra.  Since there is also an analogous root order, it is
natural to hope that a similar bijection can be used to understand
these regions.  We show that such a bijection does hold for type
$H_3$ and for type $I_2(m)$, including arbitrary ratio of root
lengths when $m$ is even, but does \textit{not}\/ hold for type
$H_4$.  We give a criterion that explains this failure and
a list of the 16 antichains in the $H_4$ root
order which correspond to empty regions.
\end{abstract}

\begin{center}
\end{center}


\section{Introduction}

Dominant regions in the Shi or Catalan hyperplane arrangement
associated to a crystallographic root system, $\Phi$, have
appeared in several contexts.  In his work on Kazhdan-Lusztig
cells for affine Weyl groups~\cite{Sh1, Sh2}, Jian-Yi Shi used
sign types to index and count these regions.  The sign type
encodes where the region lies relative to the affine hyperplanes.
In later combinatorial work, see~\cite[Remark 2]{Rei}, this result
was reformulated as a bijection between dominant regions and
antichains in the poset of positive roots in $\Phi$ under the root
order.  More recently, a uniform (case-free) proof of the
bijection was given in the context of work on ad-nilpotent ideals
in Lie algebras~\cite{CP}.  Uniform proofs have also been given
for a similar bijection and counting formula for bounded dominant
regions, as well as for various refined counts~\cite{A2, P, S}.

Representative points in the dominant regions are also used in
obtaining the support of spherical unitary representations of real
and $p$-adic algebraic groups.  This is accomplished by proving
that it suffices to test unitarity of representations of
associated affine and graded Hecke algebras~\cite{BM1, BM2} and
performing the required calculations in those settings~\cite{BM3,
Ba1, Ba2, C, V}.

When $\Phi$ is noncrystallographic, there is no corresponding
algebraic group or affine Hecke algebra, but there does exist an
associated graded Hecke algebra and the same affine hyperplanes
arise naturally as the support of its reducible principal series
modules~\cite{KR}.  Analogous calculations, involving dominant
regions in the same type of hyperplane arrangement, should be
involved in finding the support of spherical unitary
representations of these algebras. Since there is also a natural
extension of the root order to the noncrystallographic setting, a
first step is to investigate whether antichains in this root order
are still in bijection with dominant regions. We do so here,
motivated also by recent interest in generalizing various
combinatorial objects counted by Catalan numbers to
noncrystallographic type~\cite{AIM}.

We begin with a few preliminaries in order to state the main
results of the paper (for greater detail see
Section~\ref{se:Properties}).  Let $\Phi$ be a (not necessarily
crystallographic) root system in a Euclidean vector space $V$ and
fix a set of positive roots, $\Phi^+$, and simple roots,
$\Delta=\{\alpha_1,\alpha_2,\dots,\alpha_n\}$. Define the
hyperplanes $H_{\alpha,c}=\{v\in V \mid (\alpha \mid v)=c\}$. Then
the connected components of
\[V\setminus \bigcup\limits_{\substack{\beta\in\Phi^+ \\ c=-1,0,1}}
H_{\beta,c}\] contained in the fundamental chamber are
called \textit{dominant regions}. If $\Phi$ is crystallographic
this union of hyperplanes is called the Catalan hyperplane arrangement.

Define a partial order on $\Phi^+$, which we will call the
\textit{root order}, by
$$\beta \leq \gamma \qquad\hbox{ if }
\gamma-\beta=\sum\limits_{i=1}^n c_i \alpha_i \hbox{ with } c_i
\in \mathbb{R}_{\geq 0}.$$  Each increasing set in the
poset $(\Phi^+,\leq)$ is generated by its minimal elements, which
form an \textit{antichain} or incomparable set, so increasing sets
and antichains are naturally in bijection. Note that if $\Phi$ is
crystallographic, then the $c_i$ will lie in $\mathbb{Z}_{\geq 0}$.
Further properties of increasing sets and antichains are discussed in
Section~\ref{se:Properties}.

Antichains provide a means to generalize
\textit{nonnesting partitions} from symmetric groups (type $A$) to
general Weyl groups, as is attributed to Postnikov in~\cite{Rei}.
In particular, a nonnesting
partition of the set $\{1,2,\dots,n\}$ is the same as an
antichain in the root order on
$$\Phi^+_{A_{n-1}}=\{\epsilon_i-\epsilon_j \mid 1 \leq i<j \leq n\},$$
where $\{\epsilon_1,\epsilon_2,\dots,\epsilon_n\}$ is the standard
basis of $\mathbb{R}^n$, if the root
$\epsilon_i-\epsilon_j$ is assigned to each edge or arc from $i$ to $j$.

Reiner further states that Postnikov proved, in classical type,
that the nonnesting partitions are in bijection with dominant
regions in the Catalan arrangement by the map sending
an increasing set $I$ to the region defined by
\[ R_I:=\left\{\, v \in V \left|
\begin{array}{rl}
 (v \mid \alpha)>1   &\mbox{ for all } \alpha \in I \,\mbox{ and }\, \\
 0<(v \mid \alpha)<1 &\mbox{ for all }
      \alpha \in \Phi^+ \smallsetminus I \,
\end{array}\right.\right\}\,.  \]
The ideas underlying this bijection first appear in~\cite{Sh2}.

In Section~\ref{se:Non-emptyRegions}, Theorem~\ref{thm:main} we
prove that the analogous bijection holds in type $H_3$ and for all
dihedral types, $I_2(m)$.  When $m$ is even there are two orbits
of roots, and the bijection continues to hold, even as the length
of roots in one orbit varies (see Figure~\ref{fi:I2(6)} in
Section~\ref{se:Non-emptyRegions}.)

The statement and proof use both the partial order by inclusion of
antichains, denoted $\subseteq$, and (implicitly) the partial
order on antichains induced by inclusion of their corresponding
increasing sets, denoted $\preceq$.
Corollary~\ref{cor:reduction-thru-antichain} provides a sufficient
condition, given a $\subseteq$-maximal antichain, $\Lambda$, for
certain $\preceq$-subantichains of $\Lambda$ to correspond to
nonempty regions. For $H_3$ and $I_2(m)$, it is then easy to check
that all $\subseteq$-maximal antichains satisfy the condition, and
to use Corollary~\ref{cor:reduction-thru-antichain} to conclude
that all antichains correspond to nonempty regions.

However, as shown in Section~\ref{se:H4}, for type $H_4$ there are
$\subseteq$-maximal antichains that do not satisfy the condition
in Corollary~\ref{cor:reduction-thru-antichain} and the bijection
does \textit{not}\/ hold. While every dominant region gives rise
to an antichain, not all of the 429 antichains define a nonempty
region.  We use Corollary~\ref{cor:reduction-thru-antichain} to
determine easily that 401 of the antichains do correspond to
nonempty regions.  Calculations described in the tables in the
Appendix then show that exactly 16 of the remaining 28 antichains
yield empty regions.  These are listed in Theorem~\ref{thm:H4}.

Theorem~\ref{thm:criterion} in Section~\ref{se:Criterion} provides
an elementary criterion for the bijection to hold that applies to
all root systems with arbitrary root lengths. We prove there is a
bijection between dominant regions and antichains in $\Phi^+$ if
and only if for every antichain $\Lambda\neq \emptyset$, the
system $(v,\beta)=1$ for $\beta \in \Lambda$ has a solution in the
fundamental chamber.  This explains the failure of the bijection
in type $H_4$. Theorem~\ref{thm:criterion} also gives another
proof of the bijection for crystallographic root systems and
provides a way to determine if the bijection holds for nonstandard
root lengths.

These results have two important consequences.  Primarily, they
will aid in calculation of the support of the spherical unitary
dual of noncrystallographic graded Hecke algebras. Particularly
useful is the combinatorial description of the continuously
varying geometry in the case when $m$ is even, as illustrated in
Figure~\ref{fi:I2(6)} in Section~\ref{se:Non-emptyRegions}.
Secondly, they indicate one apparent obstacle to generalizing
nonnesting partitions to noncrystallographic cases.

We elaborate further on the second consequence.  If $\Phi$ is
crystallographic, there are other interesting objects that appear
likely to be related to nonnesting partitions.  For example, if
$\Phi$ has Weyl group $W$ with rank $n$, Coxeter number $h$, and
exponents $e_1, e_2, \dots, e_n$, then the numbers of associated
nonnesting partitions~\cite{A1}, noncrossing partitions~\cite{M},
vertices in simplicial associahedra and clusters~\cite{FZ} are all
given by the generalized Catalan number
$$\mathrm{Cat}(W)=\prod\limits_{i=1}^n \frac{h+e_i+1}{e_i+1}.$$
When $W=S_n$ this formula reduces to the standard Catalan number,
$\frac{1}{n+1} \binom{2n}{n}$, and the fact that this counts
noncrossing partitions goes back to Kreweras~\cite{K}.

The poset of noncrossing partitions has been generalized to
arbitrary finite Coxeter groups in independent work of David
Bessis~\cite{Be} and Tom Brady~\cite{Br, BW1}.  The
original definition of the simplicial associahedron due to Fomin
and Zelevinsky~\cite{FZ} also makes sense for any finite Coxeter
group (see~\cite[Section 5.3]{FR1}). For a noncrystallographic
reflection group $W$, the number of noncrossing partitions and the
number of vertices in associahedra still agree and are given by
$\mathrm{Cat}(W)$. Precise connections between all of these
objects are not well understood, even in the setting of Weyl
groups, though there has been some recent progress in this
direction~\cite{BW2, Ch1, Ch2, FR2, Rea}.

One question of interest is to similarly generalize nonnesting
partitions or to better understand the obstacles to doing
so~\cite[Problem 1.2]{AIM}. Our results show that antichains in
the partial order on the roots that is used here do \textit{not}
provide such a generalization. The numbers of antichains or
dominant regions are never generalized Catalan numbers except in
type $I_2(m)$ when $m=2,3,4,6$ and the root length ratio yields a
crystallographic root system (compare Table~\ref{ta:numbers} and
Table~\ref{ta:Catalan}).  We hope that these numbers and the
criterion for the failure of the bijection in type $H_4$ will be
of some help in formulating the desired generalization.

\medno\textbf{Acknowledgements.} We thank Vic Reiner for helpful
discussions and the referees for suggestions that helped improve
the exposition.  The second author thanks the American Institute
of Mathematics for support to attend the workshop ``Braid Groups,
Clusters, and Free Probability''.

\section{Properties of Increasing Sets and Antichains}
\label{se:Properties}

Let $\Phi$ be a \emph{root system} in an $n$-dimensional Euclidean
space $V$, endowed with an inner product $(\cdot \mid \cdot)$ as
defined in~\cite{Humphreys:Reflection Groups}.  In particular note
that $\Phi$ need not be crystallographic.  Let
 $\Delta=\{ \alpha_1,\alpha_2,\ldots,\alpha_n \}$ be a choice of
\emph{simple roots} with resulting \emph{positive roots}
\[ \textstyle{
   \Phi^+=\{ \beta \in \Phi
   \mid \beta=\sum_{i=1}^n c_i \alpha_i \,
   \mbox{ with } c_i \in \mathbb{R}_{\geq 0} \} }. \]
The dual basis to $\Delta$ will be denoted by $\omega_1,
\omega_2, \ldots, \omega_n$ and referred to as \emph{fundamental
weights}.

\subsection{Dominant regions}
Let $\beta \in \Phi^+$ and $c=-1,0,1$. Define the central ($c=0$)
and affine ($c=\pm 1$) hyperplanes:
\[ H_{\beta,c}=\{ v \in V \mid (v \mid \beta)=c \}. \]
A \emph{region} in $V$ is a connected component of the set
$V\setminus \bigcup\limits_{\substack{\beta\in\Phi^+ \\ c=-1,0,1}}
H_{\beta,c}$
under the Euclidean topology,
and a region is said to be \emph{dominant} if it lies in the
\emph{dominant chamber},
\begin{align*}
C&=\{ v \in V \mid (v \mid \alpha_i)>0
              \mbox{ for } 1 \leq i \leq n \} \\
 &=\{ \textstyle{\sum_{i=1}^n x_i \omega_i}
      \mid x_i>0 \, \mbox{ for } 1 \leq i \leq n \}.
\end{align*}
Define the positive side and the negative side of $H_{\beta,1}$ by
\[  H_{\beta,1}^+=\{ v \in V \mid (v \mid \beta)>1 \}
    \quad \mbox{ and } \quad
    H_{\beta,1}^-=\{ v \in V \mid (v \mid \beta)<1 \}. \]
Then every dominant region $R$ can be uniquely expressed as
\[ R=C \cap \left(
     \textstyle{\bigcap_{\beta \in \Phi^+} H_{\beta,1}^{n_\beta}}
     \right), \quad \hbox{where $n_\beta$ represents $+$ or $-$.}\]

Although consideration of the dominant regions does not require
use of the $H_{\beta,-1}$, these are included for symmetry and
because they arise in the representation theory of
noncrystallographic graded Hecke algebras.

\subsection{Increasing sets}
Define a partial order on $\Phi^+$, which we will call the
\textit{root order}, by
$$\beta \leq \gamma \qquad\hbox{ if }
\gamma-\beta=\sum\limits_{i=1}^n c_i \alpha_i \hbox{ with } c_i
\in \mathbb{R}_{\geq 0}.$$  There exist other possible orders on
$\Phi^+$ but we refer to $(\Phi^+,\leq)$ as the
\emph{positive root poset}.
An \emph{increasing set} in $(\Phi^+,\leq)$ is a subset $I$ of $\Phi^+$
such that if $\beta,\,\gamma\in \Phi^+$ satisfy
 $\beta \leq \gamma$ and $\beta \in I$, then $\gamma \in I$.
For crystallographic root systems these sets describe root spaces
of ad-nilpotent ideals, so are called ideals, but in general
posets they are called dual order ideals.

\subsection{Relating dominant regions and increasing sets}
Let $\mathcal{R}$ be the set of all dominant regions in $V$.
 To each $R \in \mathcal{R}$, assign the subset of $\Phi^+$:
\[ I_R=\{ \beta \in \Phi^+ \mid R \subseteq H_{\beta,1}^+ \}. \]
 Let $\mathcal{I}$ be the set of all increasing sets in $(\Phi^+,\leq)$.
 To each $I \in \mathcal{I}$, assign the subset of $C$:
\[ R_I=\{ v \in C
   \mid (v \mid \beta)>1 \mbox{ for all } \beta \in I \,
   \mbox{ and } \,
        (v \mid \gamma)<1 \mbox{ for all } \gamma \in I^c \}, \]
where $I^c$ represents the complement set of $I$ in $\Phi^+$,
i.e., $I^c=\Phi^+ \smallsetminus I$.

It is easy to show that for $\beta,\, \gamma \in \Phi^+$,
\begin{equation} \tag{*} \label{eq:compare-inner-product}
 \hbox{if  $\beta \leq \gamma$
 (resp. $\beta < \gamma$),
 then $(v \mid \beta) \leq (v \mid \gamma)$
 (resp. $(v \mid \beta)<(v \mid \gamma)$)
 for $v \in C$.}
\end{equation}

\begin{lemma}  \label{lem:region-to-increasing set}
If $R \in \mathcal{R}$, then $I_R$ is an increasing set in $(\Phi^+,\leq)$
and $R_{I_R}=R$.
\end{lemma}

\begin{proof}
Let $\beta,\,\gamma \in \Phi^+$ be such that
 $\beta \leq \gamma$ and $\beta \in I_R$.
Then $R \subseteq H_{\beta,1}^+$, or equivalently,
 $(v \mid \beta)>1$ for any $v \in R$. Note that $R \subseteq C$
and $\beta \leq \gamma$. By~(\ref{eq:compare-inner-product}), $\beta
\leq \gamma$ implies
 $(v \mid \gamma) \geq (v \mid \beta)>1$ for all $v \in R$,
i.e., $R \subseteq H_{\gamma,1}^+$ and $\gamma \in I_R$. So $I_R$
is an increasing set.

Let $v \in C$. Then
\begin{align*}
       v \in R
 &\iff v \in H_{\beta,1}^+ \mbox{ for all } \beta \in I_R
       \, \mbox{ and } \,
       v \in H_{\gamma,1}^- \mbox{ for all } \gamma \in (I_R)^c \\
 &\iff (v \mid \beta)>1 \mbox{ for all } \beta \in I_R
       \, \mbox{ and } \,
       (v \mid \gamma)<1 \mbox{ for all } \gamma \in (I_R)^c \\
 &\iff  v \in R_{I_R}.
\end{align*}
So $R_{I_R}=R$.
\end{proof}

\subsection{Remark}
Let the function
 $f: \mathcal{R} \rightarrow \mathcal{I}$ be defined by $f(R)=I_R$
for $R \in \mathcal{R}$. If $R_1,\,R_2 \in \mathcal{R}$ with
$f(R_1)=f(R_2)$, then $I_{R_1}=I_{R_2}$ implies $R_1=R_2$ by
Lemma~\ref{lem:region-to-increasing set}. So $f$ is always
injective. We prove $f$ is also surjective for type $H_3$ and
dihedral type $I_2(m)$ $(m \geq 2)$ in
Section~\ref{se:Non-emptyRegions} and find the image of $f$ as a
proper subset of $\mathcal{I}$ for type $H_4$ in Section~\ref{se:H4}.

\subsection{Antichains} \label{def:antichain}
An \emph{antichain} in a poset $(P,\leq)$ is a subset $\Lambda$ of
$P$ such that any distinct
 $\beta,\,\gamma \in \Lambda$
are incomparable with respect to $\leq$. Any subset of an
antichain is an antichain.

If $I$ is an increasing set in $(\Phi^+,\leq)$, then
\begin{enumerate}
\item  $I_{\min}$ represents the set of all $\leq$-minimal
elements in $I$;

\item  $I_{\max}^c$ represents the set of all $\leq$-maximal
elements in $I^c$.
\end{enumerate}
It is evident that both $I_{\min}$ and $I_{\max}^c$ are antichains
in $(\Phi^+,\leq)$. By (*), we get
\[ R_I=\{ v \in C
   \mid (v \mid \beta)>1 \mbox{ for all } \beta \in I_{\min}
   \, \mbox{ and } \,
        (v \mid \gamma)<1 \mbox{ for all } \gamma \in I_{\max}^c \}. \]

\begin{lemma} \label{lem:antichain-to-increasing set}
Let $\Lambda$ be an antichain in $(\Phi^+,\leq)$. Then
\[ I(\Lambda)=\{ \beta \in \Phi^+ \mid \delta \leq \beta \,
                 \mbox{ for some } \delta \in \Lambda \} \]
is an increasing set in $(\Phi^+,\leq)$ and
 $I(\Lambda)_{\min}=\Lambda$.
\end{lemma}

\begin{proof}
This follows easily from the definitions of increasing set and
antichain.
\end{proof}

For any antichains $\Lambda'$ and $\Lambda$ in $(\Phi^+,\leq)$, we
write $\Lambda' \preceq \Lambda$ if
 $I(\Lambda') \subseteq I(\Lambda)$.
By Lemma~\ref{lem:antichain-to-increasing set}, the binary relation
$\preceq$ is a partial order on the set of all antichains in
$(\Phi^+,\leq)$.

\begin{lemma} \label{lem:big-small-ideals}
Let $I$ be an increasing set in $(\Phi^+,\leq)$.
\begin{enumerate}
\item[(i)]  If $\Lambda \subseteq I_{\min}$, then
 $J=I \smallsetminus \Lambda$
is an increasing set in $(\Phi^+,\leq)$ and
 $\Lambda \subseteq J_{\max}^c$.

\item[(ii)]  If $\Lambda \subseteq I_{\max}^c$, then
 $J=I \cup \Lambda$
is an increasing set in $(\Phi^+,\leq)$ and
 $\Lambda \subseteq J_{\min}$.
\end{enumerate}
\end{lemma}

\begin{proof}
(i) and (ii) follow easily from the definition of increasing set.
\end{proof}

\begin{corollary} \label{cor:criterion-for-I^c_max}
If $I$ is an increasing set in $(\Phi^+,\leq)$, then
\[ I_{\max}^c=\{\, \beta \in I^c \mid I \cup \{\beta\}
      \mbox{ is an increasing set} \,\}. \]
\end{corollary}

By observing the positive root posets, e.g.
Figure~\ref{fi:H4rootposet} in Section~\ref{se:H4} and
Figure~\ref{fi:I2(6)} in Section~\ref{se:Non-emptyRegions}, it is
easy to see

\begin{lemma} \label{lem:decomposition-of-positive-roots}
The set of positive roots can be decomposed into a disjoint union
of subsets $\Phi_{\alpha_i}^+$ for $1 \leq i \leq n$, where all
the roots in each $\Phi_{\alpha_i}^+$ can be arranged in a certain
order
 $\beta_0, \beta_1, \ldots, \beta_m$
such that $\beta_0 \in \Delta$, $\beta_j=s_{\alpha} (\beta_{j-1})$
for some simple root $\alpha \in \Delta$, and $\beta_{j-1} \leq
\beta_j$ for $1 \leq j \leq m$. In particular, all the roots in
the same $\Phi^+_{\alpha_i}$ are comparable.
\end{lemma}

As an immediate consequence of
Lemma~\ref{lem:decomposition-of-positive-roots}, we have

\begin{corollary} \label{cor:size-of-antichain}
If $\Phi$ is a root system of rank $n$, then any antichain in the
positive root poset $(\Phi^+,\leq)$ has at most $n$ positive
roots.
\end{corollary}

If $\Phi$ is crystallographic there is exactly one such
antichain~\cite{P}, but if $\Phi$ is noncrystallographic there can
be several (see Figure~\ref{fi:I2(6)} in
Section~\ref{se:Non-emptyRegions}). Suppose that $V_0$ is a
subspace of the Euclidean space $V$ such that $\Phi_0=\Phi \cap
V_0$ is nonempty. Then $\Phi_0$ is a root system in $V_0$ and
there exists a set $\Delta_0$ of simple roots so that
$\Phi_0^+=\Phi^+ \cap V_0$. Let $\leq_0$ be the root order on
$\Phi_0^+$.
 If $\beta,\,\gamma \in \Phi_0^+$ satisfy $\beta \leq_0 \gamma$,
then $\beta \leq \gamma$.

\begin{proposition}  \label{prop:linear-independence-of-antichain}
If $\Lambda$ is an antichain in $(\Phi^+,\leq)$, then $\Lambda$ is
a linearly independent set.
\end{proposition}

\begin{proof}
Let $n$ be the rank of $\Phi$. If $n=1$, then $|\Lambda| \leq 1$
and hence $\Lambda$ is linearly independent. Suppose that for any
root system $\Phi_0$ of rank less than $n$, every antichain
$\Lambda_0$ in $(\Phi_0^+,\leq_0)$ is linearly independent. If
 $\Lambda=\{\beta_1, \ldots, \beta_{m}, \beta_{m+1}\}$
is linearly dependent, then we may assume
 $\beta_{m+1}=\sum_{i=1}^m c_i \beta_i$ for $c_i \in \mathbb{R}$.
Let $V_0$ be the linear subspace of $V$ spanned by
 $\{ \beta_1, \ldots, \beta_m \}$.
Set $\Phi_0=\Phi \cap V_0$ and $\Phi_0^+=\Phi^+ \cap V_0$.
 $\Lambda_0=\{\beta_1,\ldots,\beta_m,\beta_{m+1}\}$
is an antichain in $(\Phi_0^+,\leq_0)$. By
Corollary~\ref{cor:size-of-antichain}, we get
 $m+1 \leq n$ and $\rank(\Phi_0) \leq m <n$.
By the induction hypothesis, the antichain
 $\Lambda=\{\beta_1, \ldots, \beta_m, \beta_{m+1}\}$
is linearly independent.
\end{proof}

For crystallographic root systems, the results in
Corollary~\ref{cor:size-of-antichain} and
Proposition~\ref{prop:linear-independence-of-antichain} are proven
in~\cite[Proposition 2.10]{P} by other means.

\subsection{Definition} \label{def:Int}
Let $\Lambda$ be a nonempty antichain in $(\Phi^+,\leq)$. Define
\begin{enumerate}
\item  \hskip0.23truecm
 $\Int(\Lambda)=\{ v \in V \mid (v \mid \beta)=1
  \mbox{ for all } \beta \in \Lambda \}
  =\bigcap_{\beta \in \Lambda} H_{\beta,1}$,

\item
 $\Int_{C}(\Lambda)=\{ v \in C \mid (v \mid \beta)=1
  \mbox{ for all } \beta \in \Lambda \}
  =\Int(\Lambda) \cap C$.
\end{enumerate}

By Proposition~\ref{prop:linear-independence-of-antichain},
$\Int(\Lambda)$ is a nonempty convex set, and so is
$\Int_{C}(\Lambda)$.

\begin{proposition} \label{prop:criterion-for-incomparable}
Let $\Lambda=\{\beta,\,\gamma\}$ for some distinct
 $\beta,\,\gamma \in \Phi^+$.
Then $\Lambda$ is an antichain in $(\Phi^+,\leq)$ if and only if
$\Int_{C}(\Lambda) \ne \emptyset$.
\end{proposition}

\begin{proof}
Let
 $\beta=\sum_{i=1}^n c_i \alpha_i$ and
 $\gamma=\sum_{i=1}^n d_i \alpha_i$,
where $c_i \geq 0$ and $d_i \geq 0$.

$\Rightarrow)$  Suppose that $\Lambda$ is an antichain. Then
$\beta$ and $\gamma$ are incomparable, and there are $j,k$ such
that $c_j>d_j$ and $c_k<d_k$. Choose positive
 $x_1, x_2, \ldots, x_n \in \mathbb{R}$
such that
 $\sum_{i=1}^n x_i (c_i-d_i)=0$.
Set $v=\sum_{i=1}^n x_i\omega_i \in C$. Then
 $\frac{v}{(v \mid \beta)} \in \Int_{C}(\Lambda)$.

$\Leftarrow)$ Suppose that $v \in \Int_{C}(\Lambda)$. Let
$v=\sum_{i=1}^n x_i \omega_i$ for $x_i>0$. Then
 $\sum_{i=1}^n x_i(c_i-d_i)=0$.
There are $j,k$ such that $c_j<d_j$ and $c_k>d_k$. $\beta, \gamma$
are incomparable, so $\Lambda$ is an antichain.
\end{proof}

\section{Non-Empty Regions for Types $H_3$ and $I_2(m)$}
\label{se:Non-emptyRegions}

In the remainder of the paper, we write $B_{\epsilon}(v_0)$ for
the open ball in the Euclidean space $V$ with center $v_0 \in V$
and radius $\epsilon>0$, i.e.,
\[ B_\epsilon(v_0)
   =\{ v \in V \mid \| v-v_0 \|<\epsilon \}\,. \]

Note that $R_\emptyset$ and $R_{\Phi^{+}}$ are a nonempty dominant
regions.

\begin{theorem}  \label{thm:reduction-thru-antichain}
Let $I$ be a nonempty increasing set in $(\Phi^+,\leq)$.
\begin{enumerate}
\item[(i)]  If there exists $v_0 \in \Int_{C}(I_{\min})$ such that
$(v_0 \mid \beta)<1$ for all $\beta \in I_{\max}^c$, then
 $R_{I \smallsetminus \Lambda} \ne \emptyset$
for all $\Lambda \subseteq I_{\min}$.

\item[(ii)]  If there exists $v_0 \in \Int_{C}(I_{\max}^c)$ such that
$(v_0 \mid \beta)>1$ for all $\beta \in I_{\min}$, then
 $R_{I \cup \Lambda} \ne \emptyset$
for all $\Lambda \subseteq I_{\max}^c$.
\end{enumerate}
\end{theorem}

\begin{proof}
We prove part (i) in three steps.

(1) Let $\gamma \in \Phi^+ \smallsetminus I_{\min}$.

 If $\gamma \in I$, then there exists $\beta \in I_{\min}$
with $\beta<\gamma$, so
 $(v_0 \mid \gamma)>(v_0 \mid \beta)=1$.

 If $\gamma \in I^c$, then there exists $\beta \in I_{\max}^c$
with $\gamma \leq \beta$, so
 $(v_0 \mid \gamma) \leq (v_0 \mid \beta)<1$.

(2) Set
 $d_{\gamma}=\min \{ \| v-v_0 \| \mid v \in H_{\gamma,1} \}$
for $\gamma \in \Phi^+ \smallsetminus I_{\min}$ and
 $d_i=\min\{ \|v-v_0\| \mid v \in H_{\alpha_i,0} \}$
for $1 \leq i \leq n$. By step (1) and $v_0 \in C$,
 $\epsilon=\min\{\, d_\gamma,\,d_i \mid
  \gamma \in \Phi^+ \smallsetminus I_{\min}
  \mbox{ and } 1 \leq i \leq n \,\}>0$.
So $B_{\epsilon}(v_0)$ is contained in $C$ and does not intersect
$H_{\gamma,1}$ for $\gamma \in \Phi^+ \smallsetminus I_{\min}$. By
step (1), for any $v \in B_{\epsilon}(v_0)$ we have
 $(v \mid \gamma)>1$ if $\gamma \in I \smallsetminus I_{\min}$ and
 $(v \mid \gamma)<1$ if $\gamma \in \Phi^+ \smallsetminus I$.

(3) By Proposition~\ref{prop:linear-independence-of-antichain},
$I_{\min}$ is a linearly independent set in $V$. There exists
 $u \in V$
such that
 $(u \mid \beta)=-1 \, \mbox{ for all } \beta \in \Lambda$
and
 $(u \mid \beta)=1 \, \mbox{ for all } \beta \in
 I_{\min} \smallsetminus \Lambda$.
Set $v=v_0+\frac{\epsilon}{2} \cdot \frac{u}{\|u\|}$. Then
 $\|v-v_0\|=\frac{\epsilon}{2}<\epsilon$.
So $v \in B_{\epsilon}(v_0)$.
\begin{enumerate}
\item[]
 If $\gamma \in \Lambda$, then
 $(v \mid \gamma)
 =(v_0 \mid \gamma)+\frac{\epsilon}{2\|u\|}(u \mid \beta)
 =1-\frac{\epsilon}{2\|u\|}<1$.

\item[]
 If $\gamma \in I_{\min} \smallsetminus \Lambda$, then
 $(v \mid \beta)
 =(v_0 \mid \beta)+\frac{\epsilon}{2\|u\|}(u \mid \beta)
 =1+\frac{\epsilon}{2\|u\|}>1$.
\end{enumerate}
Therefore, $v \in R_{I \smallsetminus \Lambda}$ and
 $R_{I \smallsetminus \Lambda} \ne \emptyset$.

Part (ii) of this theorem can be proved in the same manner.
\end{proof}

Note that any antichain of cardinality $n$ in $(\Phi^+,\leq)$ is
$\subseteq$-maximal, but not every $\subseteq$-maximal antichain
in $(\Phi^+,\leq)$ must have cardinality $n$.

\begin{corollary} \label{cor:reduction-thru-antichain}
Let $I$ be a nonempty increasing set in $(\Phi^+,\leq)$.
\begin{enumerate}
\item[(i)]  If $I_{\min}$ is $\subseteq$-maximal in
$(\Phi^+,\leq)$ and $\Int_{C}(I_{\min})\neq\emptyset$, then
 $R_{I \smallsetminus \Lambda} \ne \emptyset$ for all $\Lambda
\subseteq I_{\min}$.

\item[(ii)]  If $I_{\max}^c$ is $\subseteq$-maximal in
$(\Phi^+,\leq)$ and $\Int_{C}(I_{\max}^c)\neq\emptyset$, then
 $R_{I \cup \Lambda} \ne \emptyset$
for all $\Lambda \subseteq I_{\max}^c$.
\end{enumerate}
\end{corollary}

\begin{proof}
If $I_{\min}$ (resp. $I_{\max}^{c}$) is $\subseteq$-maximal in
$(\Phi^+,\leq)$, every
 $\gamma \in \Phi^+ \smallsetminus I_{\min}$
(resp.
 $\gamma \in \Phi^+ \smallsetminus I_{\max}^{c}$)
is comparable with some $\beta \in I_{\min}$ (resp. $\beta \in
I_{\max}^{c}$). The results follow from
Theorem~\ref{thm:reduction-thru-antichain}.
\end{proof}

\subsection{Types $I_2(m)$ and $H_3$}
\label{se:H3}
For type $I_2(m)$, the Coxeter graph is
$\,\,$\begin{texdraw}
  \drawdim cm \linewd 0.02 \setgray 0.1
  \textref h:C v:C

  \htext(1 1){$\circ$} \move(1.06 1) \lvec(1.91 1)
  \htext(2 1){$\circ$} \htext(1.5 1.25){\tiny$m$}

  \htext(1 0.75){\tiny $1$} \htext(2 0.75){\tiny$2$}
\end{texdraw}$\,\,$.
There are $m$ positive roots and when $m$ is
even the positive root poset varies depending on the ratio of
root lengths.
For example, for type $I_2(6)$
with $\|\alpha_1\|=1$ and $\|\alpha_2\|=r$
the positive roots are

\vskip0.5truecm
\begin{center}
\begin{tabular}{c}
$\xymatrix{ *={} \save[]+/va(0) 2cm/ *{\beta_1}\ar @{-} []|
\restore \save[]+/va(30) 2cm/ *{\beta_2} \ar @{-} []| \restore
\save[]+/va(60) 2cm/ *{\beta_3} \ar @{-} []| \restore
\save[]+/va(90) 2cm/ *{\beta_4} \ar @{-} []| \restore
\save[]+/va(120) 2cm/ *{\beta_5} \ar @{-} []| \restore
\save[]+/va(150) 2cm/ *{\beta_6} \ar @{-} []| \restore}$
\end{tabular} \hskip2truecm
\begin{tabular}{l@{\hskip1truecm}l}
   $\beta_1=\alpha_1 \hskip1truecm$
 & $\beta_4=\sqrt{3}r\alpha_1+2\alpha_2$ \\
   $\beta_2=\sqrt{3}r\alpha_1+\alpha_2$
 & $\beta_5=\alpha_1+\frac{\sqrt{3}}{r}\alpha_2$ \\
   $\beta_3=2\alpha_1+\frac{\sqrt{3}}{r}\alpha_2$
 & $\beta_6=\alpha_2$.
\end{tabular}
\end{center}

\vskip0.5truecm

\begin{landscape}
\begin{figure}
\caption{Dominant regions and positive root posets in type
$I_2(6)$} \label{fi:I2(6)}
\includegraphics[trim=12.75cm 5cm 2.2cm 5cm,clip,width=4.5cm]{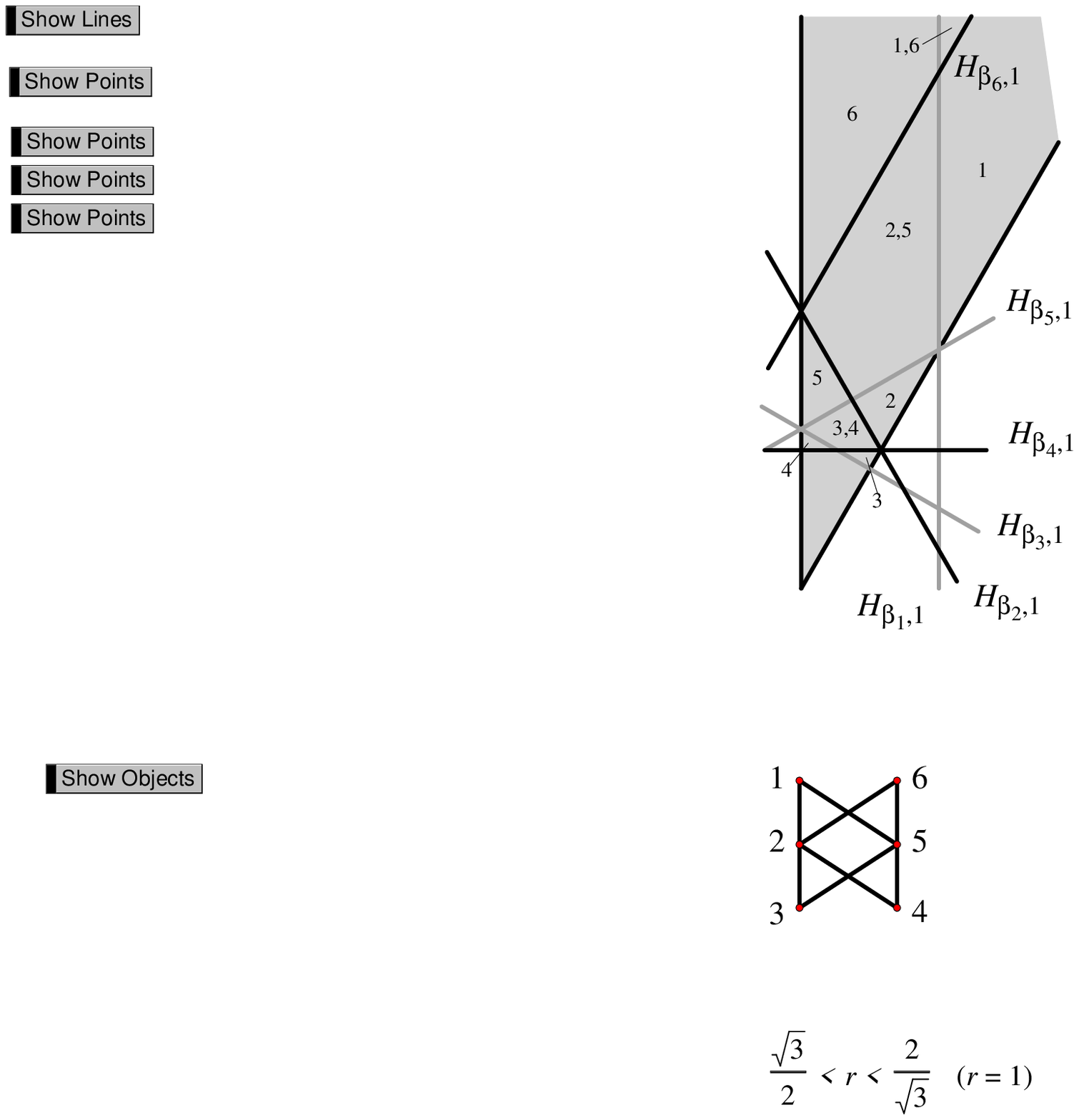}
\includegraphics[trim=12.75cm 5cm 2.2cm 5cm,clip,width=4.5cm]{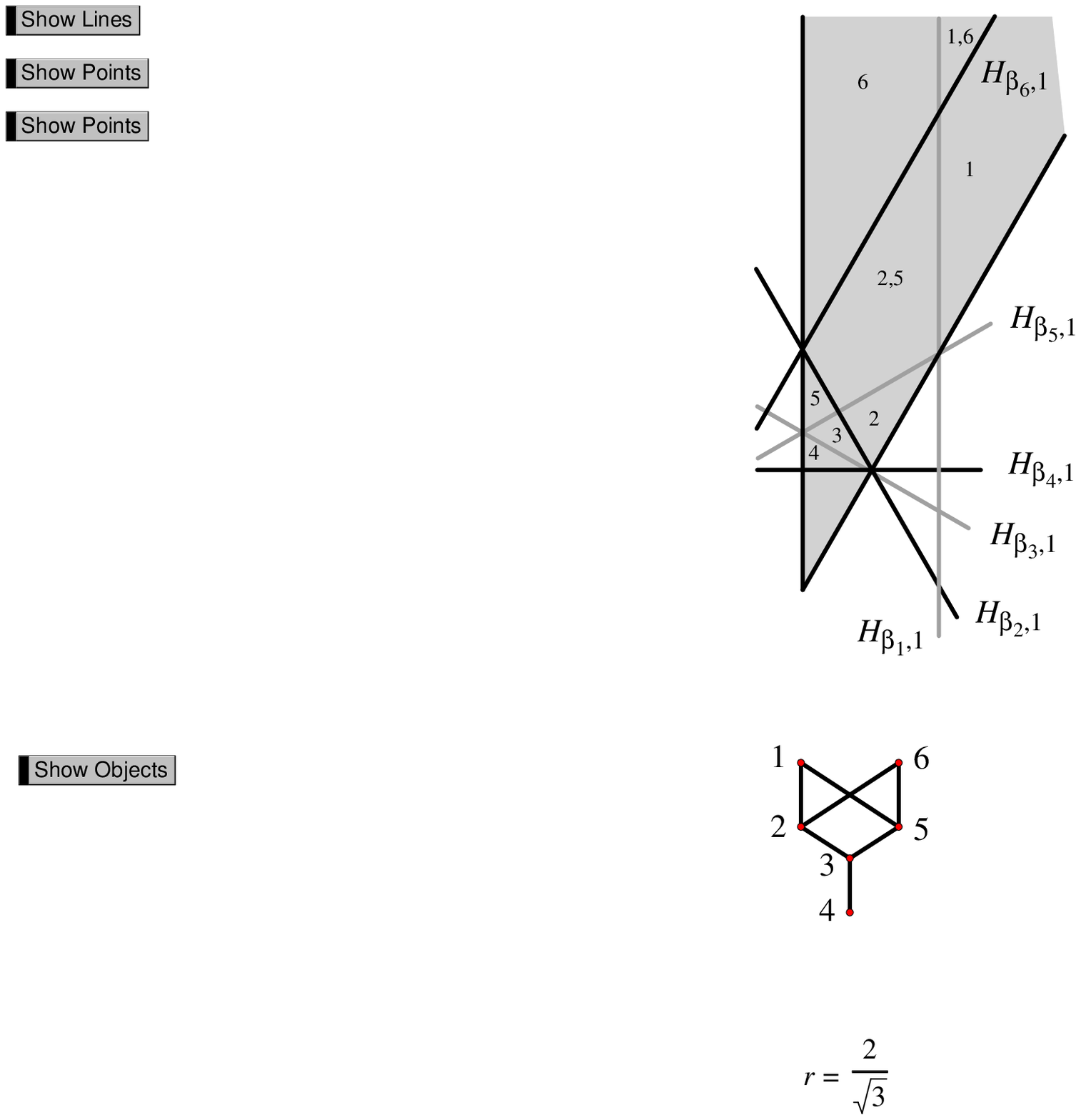}
\includegraphics[trim=12.75cm 5cm 2.2cm 5cm,clip,width=4.5cm]{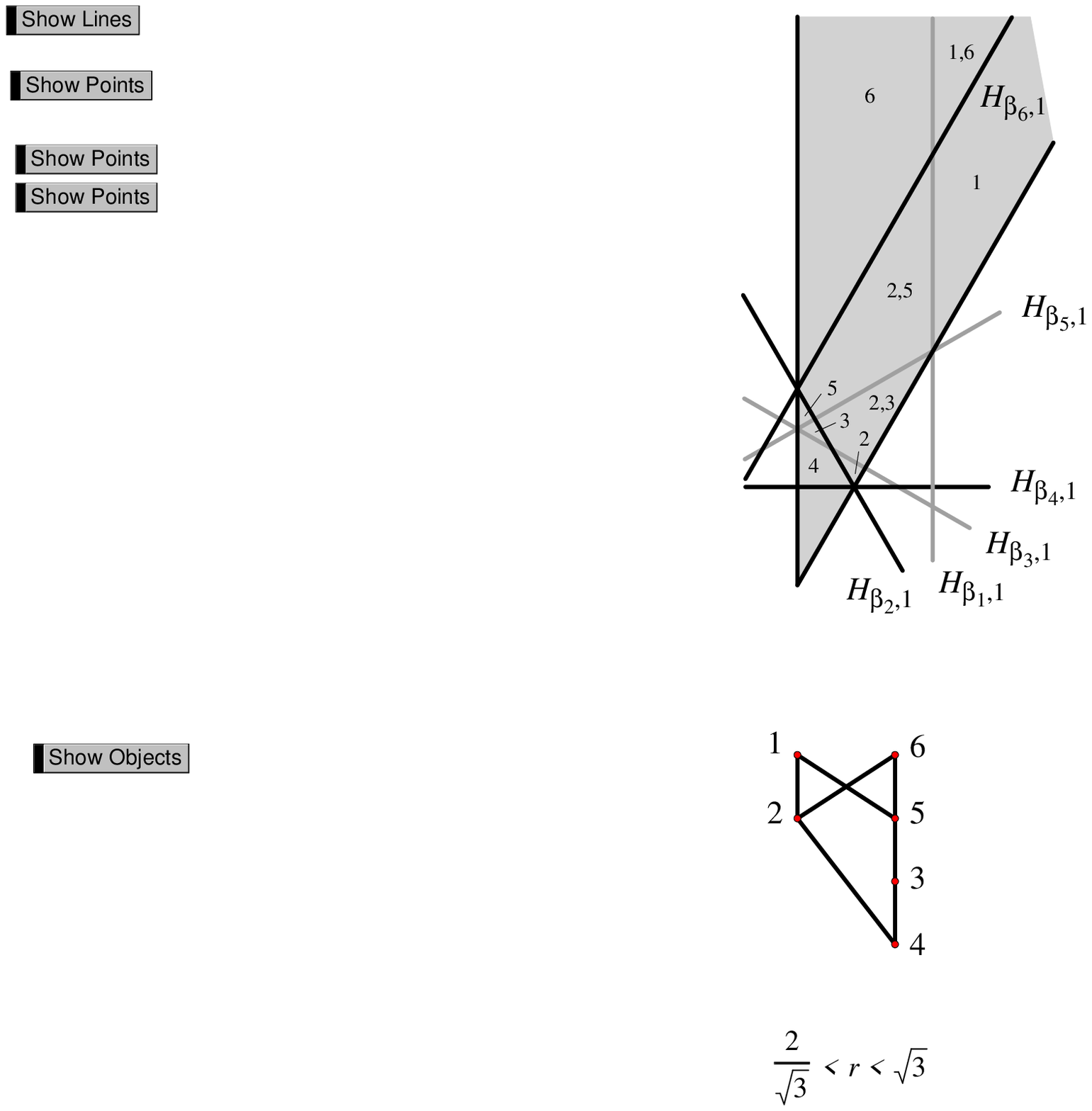}
\includegraphics[trim=12.75cm 5cm 2.2cm 5cm,clip,width=4.5cm]{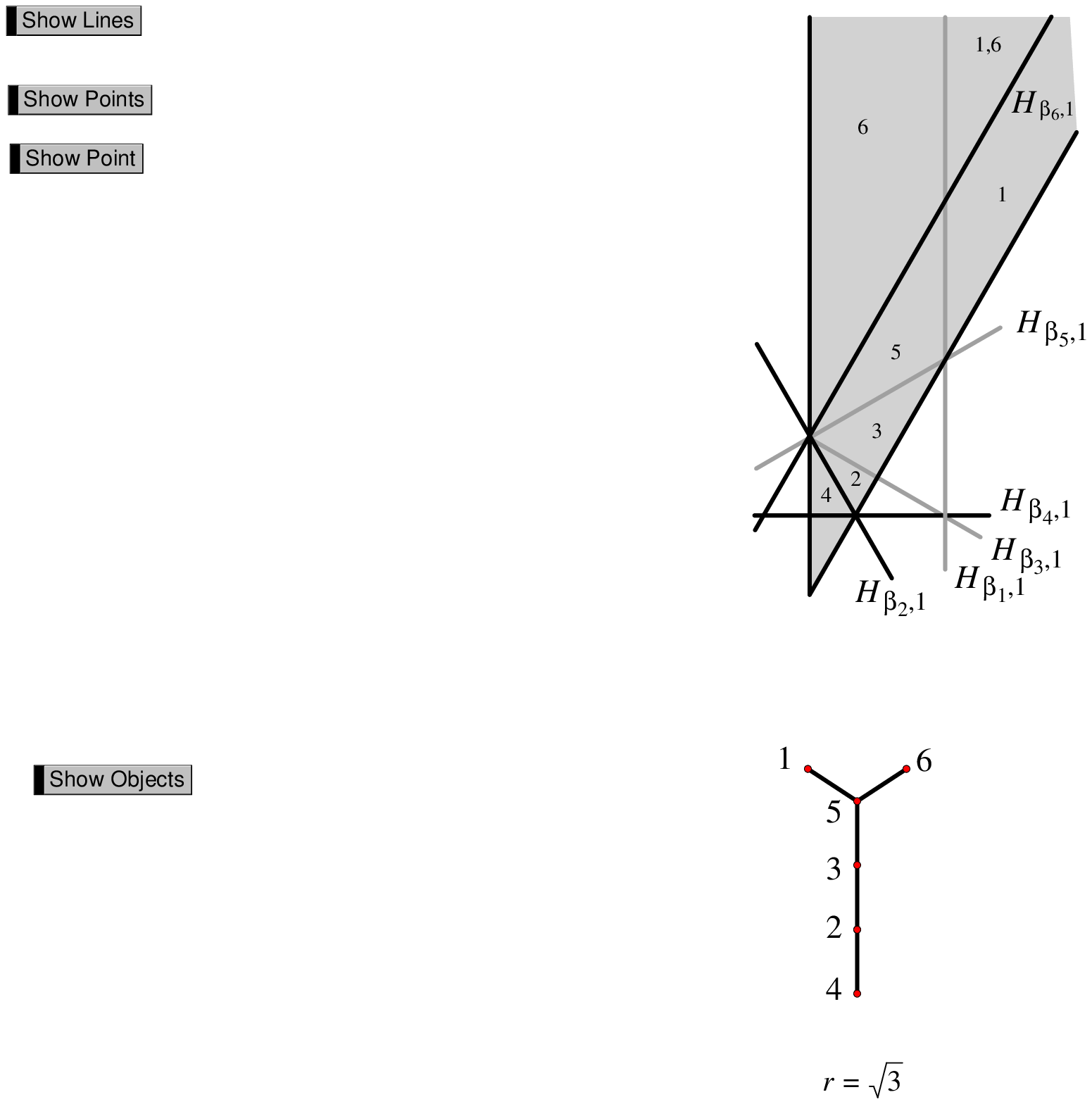}
\includegraphics[trim=12.75cm 5cm 2.2cm 5cm,clip,width=4.5cm]{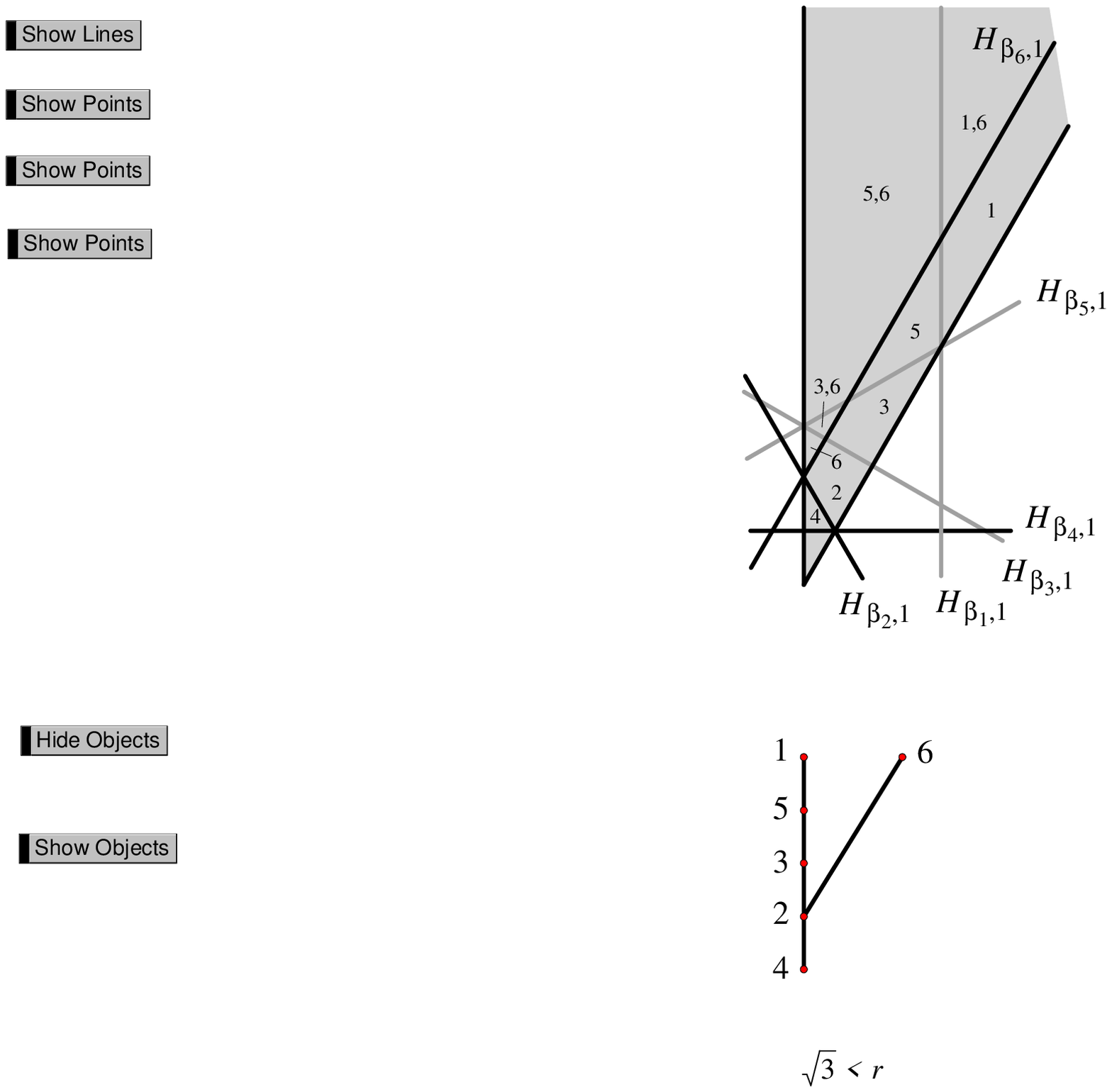}
\end{figure}
\end{landscape}

\noindent Figure~\ref{fi:I2(6)} illustrates the hyperplanes,
dominant regions and positive root posets for various values of
$r$ ($\leq$-minimal elements are at the top of the posets).  The
corresponding antichain is marked on each dominant region except
for $R_\emptyset$.  The root $\beta_i$ is abbreviated to $i$ in
both the poset and the antichain marking a region.  The four
similar but dual pictures for $r<\frac{\sqrt{3}}{2}$ are omitted.
Note that regardless of whether $m$ is even or what $r$ is,
$\subseteq$-maximal antichains are always of size at most two.

For type $H_3$, the Coxeter graph is $\,\,$
\begin{texdraw}
  \drawdim cm \linewd 0.02 \setgray 0.1
  \textref h:C v:C
  \htext(1 1){$\circ$} \move(1.06 1) \lvec(1.91 1)
  \htext(2 1){$\circ$} \htext(1.5 1.25){\tiny$5$} \move(2.06 1) \lvec(2.91 1)
  \htext(3 1){$\circ$}
  \htext(1 0.75){\tiny $1$} \htext(2 0.75){\tiny$2$} \htext(3 0.75){\tiny $3$}
\end{texdraw}
$\,\,$. There are $15$ positive roots and the positive root poset
for type $H_3$ can be obtained from $(\Phi^+,\leq)$ for type $H_4$
by restricting to the roots with zero coefficient on $\alpha_4$.
This can be seen in Figure~\ref{fi:H4rootposet}, where
\[\tau=2\cos(\pi/5)=\frac{1+\sqrt{5}}{2},\]
\noindent $\leq$-minimal elements are at the top, the values given
are the coefficients of $\alpha_1, \alpha_2, \alpha_3$, and
$\alpha_4$ respectively, and covering relations that do not arise
from simple reflections are marked with dashed lines.  The
positive roots for $H_3$ are numbered consecutively across rows
from top to bottom. Suppose that $I_{\min}$ is a three-element
antichain. We will use the notation $c_1,\,c_2,\,c_3$ to denote a
solution of the form $c_1\omega_1+c_2\omega_2+c_3\omega_3$ to the
system
 $(v \mid \beta)=1$ for all $\beta \in I_{\min}$.

\begin{center}
\begin{tabular}{|c|c|}
\hline
 $I_{\min}$ & Solution to $(v \mid \beta)=1$ for all $\beta \in I_{\min}$
\\ \hline \hline

  $\{ \alpha_1,\alpha_2,\alpha_3 \}$
 & $1, \, 1, \, 1$
\\ \hline

  $\{ \alpha_3,\alpha_4,\alpha_5 \}$
 & $2-\tau, \, 2-\tau, \, 1$
\\ \hline

  $\{ \alpha_4,\alpha_5,\alpha_6 \}$
 & $2-\tau, \, 2-\tau, \, \tau-1$
\\ \hline

  $\{ \alpha_7,\alpha_8,\alpha_9 \}$
 & $2-\tau, \, 2\tau-3, \,5-3\tau$
\\ \hline
\end{tabular}
\end{center}

\begin{theorem} \label{thm:main}
Let $\Phi$ be a root system of type $H_3$ or of dihedral type
$I_2(m)$ $(m \geq 2)$. Then the function
 $f: \mathcal{R} \rightarrow \mathcal{I}$
given by $f(R)=I_R$ for $R \in \mathcal{R}$ is a bijection.
\end{theorem}

\begin{proof}
Let $I$ be an increasing set in $(\Phi^+,\leq)$ which is not
$\emptyset$ or $\Phi^{+}$. Set
 $\Lambda=I_{\max}^c$ and $J=I \cup \Lambda$,
where $J$ is an increasing set by Lemma~\ref{lem:big-small-ideals}(ii).

We claim that $J_{\min}$ is a $\subseteq$-maximal antichain in
$(\Phi^+,\leq)$. It suffices to show that for any $\gamma \in
J^c$, there is $\beta \in J_{\min}$ such that $\gamma \leq \beta$.
 Since $\gamma \in J^c \subseteq I^c$,
there exists some $\beta \in I_{\max}^c=\Lambda$ such that
 $\gamma \leq \beta$.
We also have $\beta \in \Lambda \subseteq J_{\min}$ by
Lemma~\ref{lem:big-small-ideals}(ii).

By Proposition~\ref{prop:criterion-for-incomparable} and
Section~\ref{se:H3}, $\Int_{C}(J_{\min})\neq\emptyset$ and then by
Corollary~\ref{cor:reduction-thru-antichain}(i), we obtain
$R_I=R_{J \smallsetminus \Lambda} \ne \emptyset$. Since
  $f(R_I)=I_{R_I}
         =\{ \beta \mid R_I \subseteq H_{\beta,1}^+\}
         =I$,
$f$ is surjective.
\end{proof}

\section{Empty Regions and Nonempty Regions for Type $H_4$}
\label{se:H4}

\begin{figure}
\caption{Positive root poset for type $H_4$
(bottom is bent to fit on the page)}
\label{fi:H4rootposet}
\tiny{$$\xymatrix@1{
*+[F-:<2pt>]{1,0,0,0} \ar@{-}[d]^{s_2}\ar@{--}[dr] &
*+[F-:<2pt>]{0,1,0,0} \ar@{--}[dl]\ar@{-}[d]^{s_1}\ar@{-}[dr]^{s_3} &
*+[F-:<2pt>]{0,0,1,0} \ar@{-}[d]^{s_2}\ar@{-}[dr]^{s_4} &
*+[F-:<2pt>]{0,0,0,1} \ar@{-}[d]^{s_3} \\
*+[F-:<2pt>]{1,\tau,0,0} \ar@{-}[d]^{s_1}\ar@{-}[dr]^(.4){s_3} &
*+[F-:<2pt>]{\tau,1,0,0} \ar@{-}[dl]_(.4){s_2}\ar@{-}[dr]^{s_3} &
*+[F-:<2pt>]{0,1,1,0} \ar@{--}[dl]\ar@{-}[d]^{s_1}\ar@{-}[dr]^{s_4} &
*+[F-:<2pt>]{0,0,1,1} \ar@{-}[d]^{s_2} \\
*+[F-:<2pt>]{\tau,\tau,0,0} \ar@{-}[d]^{s_3}\ar@{--}[drr] &
*+[F-:<2pt>]{1,\tau,\tau,0} \ar@{-}[dl]^{s_1}\ar@{-}[d]^(.7){s_4} &
*+[F-:<2pt>]{\tau,1,1,0} \ar@{-}[d]^(.6){s_2}\ar@{-}[dr]_(.53){s_4}\ar@{--}[dll] &
*+[F-:<2pt>]{0,1,1,1} \ar@{--}[dll]\ar@{-}[d]^{s_1} \\
*+[F-:<2pt>]{\tau,\tau,\tau,0} \ar@{-}[d]^{s_2}\ar@{-}[dr]^{s_4} &
*+[F-:<2pt>]{1,\tau,\tau,\tau} \ar@{-}[d]^(.7){s_1} &
*+[F-:<2pt>]{\tau,\tau+1,1,0} \ar@{-}[dll]_(.24){s_3}\ar@{-}[d]^(.7){s_1}\ar@{-}[dr]_(.53){s_4} &
*+[F-:<2pt>]{\tau,1,1,1} \ar@{-}[d]^{s_2}\ar@{--}[dll] \\
*+[F-:<2pt>]{\tau,\tau+1,\tau,0} \ar@{-}[d]^{s_1}\ar@{-}[dr]^{s_4}\ar@{--}[drrr] &
*+[F-:<2pt>]{\tau,\tau,\tau,\tau} \ar@{-}[d]^(.7){s_2} &
*+[F-:<2pt>]{\tau+1,\tau+1,1,0} \ar@{-}[dll]_(.24){s_3}\ar@{-}[d]^(.3){s_4} &
*+[F-:<2pt>]{\tau,\tau+1,1,1} \ar@{-}[dl]^{s_1} \ar@{-}[d]^{s_3}\ar@{--}[dll] \\
*+[F-:<2pt>]{\tau+1,\tau+1,\tau,0} \ar@{-}[d]^{s_2} \ar@{-}[dr]^{s_4}\ar@{--}[drrr] &
*+[F-:<2pt>]{\tau,\tau+1,\tau,\tau} \ar@{-}[d]^{s_1}\ar@{-}[dr]^{s_3} &
*+[F-:<2pt>]{\tau+1,\tau+1,1,1} \ar@{-}[dr]^(.4){s_3}\ar@{--}[dl] &
*+[F-:<2pt>]{\tau,\tau+1,\tau+1,1} \ar@{-}[dl]_(.4){s_4}\ar@{-}[d]^{s_1} \\
*+[F-:<2pt>]{\tau+1,2\tau,\tau,0} \ar@{-}[d]^{s_4}\ar@{--}[drr] &
*+[F-:<2pt>]{\tau+1,\tau+1,\tau,\tau} \ar@{-}[dl]^{s_2}\ar@{-}[d]^(.7){s_3} &
*+[F-:<2pt>]{\tau,\tau+1,\tau+1,\tau} \ar@{-}[dl]_{s_1} &
*+[F-:<2pt>]{\tau+1,\tau+1,\tau+1,1} \ar@{-}[dll]_{s_4}\ar@{-}[dl]^{s_2} \\
*+[F-:<2pt>]{\tau+1,2\tau,\tau,\tau} \ar@{-}[d]^{s_3}\ar@{--}[dr] &
*+[F-:<2pt>]{\tau+1,\tau+1,\tau+1,\tau} \ar@{--}[dl]\ar@{-}[d]^{s_2} &
*+[F-:<2pt>]{\tau+1,2\tau+1,\tau+1,1} \ar@{-}[dl]^{s_4}\ar@{-}[d]^{s_1} \\
*+[F-:<2pt>]{\tau+1,2\tau,2\tau,\tau} \ar@{-}[d]^{s_2} &
*+[F-:<2pt>]{\tau+1,2\tau+1,\tau+1,\tau} \ar@{-}[dl]^{s_3}\ar@{-}[d]^{s_1} &
*+[F-:<2pt>]{2\tau+1,2\tau+1,\tau+1,1} \ar@{-}[dl]^{s_4}\ar@{-}[d]^{s_2} \\
*+[F-:<2pt>]{\tau+1,2\tau+1,2\tau,\tau} \ar@{-}[d]^{s_1} &
*+[F-:<2pt>]{2\tau+1,2\tau+1,\tau+1,\tau} \ar@{-}[dl]^{s_3}\ar@{-}[d]^{s_2} &
*+[F-:<2pt>]{2\tau+1,2\tau+2,\tau+1,1} \ar@{-}[dl]^{s_4}\ar@{-}[d]^{s_3} \\
*+[F-:<2pt>]{2\tau+1,2\tau+1,2\tau,\tau} \ar@{-}[d]^{s_2}\ar@{--}[dr]\ar@{--}[drr] &
*+[F-:<2pt>]{2\tau+1,2\tau+2,\tau+1,\tau} \ar@{--}[dl]\ar@{-}[d]^(.7){s_3}\ar@{--}[dr] &
*+[F-:<2pt>]{2\tau+1,2\tau+2,\tau+2,1} \ar@{--}[dl]\ar@{-}[d]^{s_4} \\
*+[F-:<2pt>]{2\tau+1,3\tau+1,2\tau,\tau} \ar@{-}[d]^{s_1}\ar@{-}[dr]^{s_3} &
*+[F-:<2pt>]{2\tau+1,2\tau+2,2\tau+1,\tau} \ar@{-}[d]^{s_2}\ar@{-}[dr]^{s_4} &
*+[F-:<2pt>]{2\tau+1,2\tau+2,\tau+2,\tau+1} \ar@{-}[d]^{s_3} \\
*+[F-:<2pt>]{2\tau+2,3\tau+1,2\tau,\tau} \ar@{-}[d]^{s_3} &
*+[F-:<2pt>]{2\tau+1,3\tau+1,2\tau+1,\tau} \ar@{-}[dl]^{s_1}\ar@{-}[d]^{s_4} &
*+[F-:<2pt>]{2\tau+1,2\tau+2,2\tau+1,\tau+1} \ar@{-}[dl]^{s_2} \\
*+[F-:<2pt>]{2\tau+2,3\tau+1,2\tau+1,\tau} \ar@{-}[d]^{s_2}\ar@{-}[dr]^{s_4} &
*+[F-:<2pt>]{2\tau+1,3\tau+1,2\tau+1,\tau+1} \ar@{-}[d]^{s_1} \\
*+[F-:<2pt>]{2\tau+2,3\tau+2,2\tau+1,\tau} \ar@{-}[d]^{s_1}\ar@{-}[dr]^{s_4} &
*+[F-:<2pt>]{2\tau+2,3\tau+1,2\tau+1,\tau+1} \ar@{-}[d]^{s_2} \\
*+[F-:<2pt>]{3\tau+1,3\tau+2,2\tau+1,\tau} \ar@{-}[d]^{s_4} &
*+[F-:<2pt>]{2\tau+2,3\tau+2,2\tau+1,\tau+1} \ar@{-}[dl]^{s_1}\ar@{-}[d]^{s_3} \\
*+[F-:<2pt>]{3\tau+1,3\tau+2,2\tau+1,\tau+1} \ar@{-}[d]^{s_3} &
*+[F-:<2pt>]{2\tau+2,3\tau+2,2\tau+2,\tau+1} \ar@{-}[dl]^{s_1} \\
*+[F-:<2pt>]{3\tau+1,3\tau+2,2\tau+2,\tau+1} \ar@{-}[d]^{s_2} \\
*+[F-:<2pt>]{3\tau+1,3\tau+3,2\tau+2,\tau+1} \ar@{-}[d]^{s_1} \\
*+[F-:<2pt>]{3\tau+2,3\tau+3,2\tau+2,\tau+1} \ar@{-}[r]^{s_2} &
*+[F-:<2pt>]{3\tau+2,4\tau+2,2\tau+2,\tau+1} \ar@{-}[r]^{s_3} &
*+[F-:<2pt>]{3\tau+2,4\tau+2,3\tau+1,\tau+1} \ar@{-}[r]^{s_4} &
*+[F-:<2pt>]{3\tau+2,4\tau+2,3\tau+1,2\tau} }$$}
\end{figure}
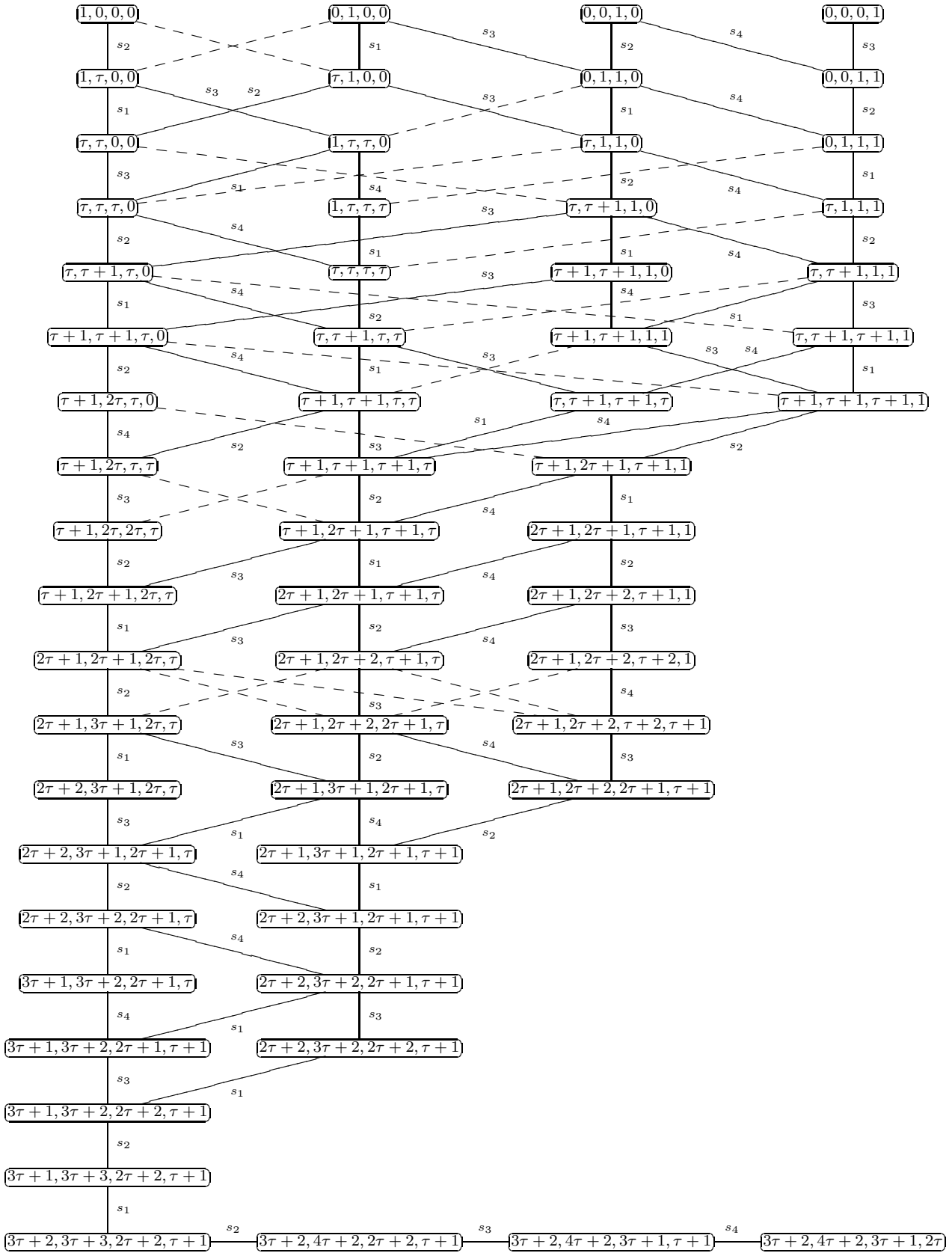

Throughout this section assume that $\Phi$ is a root system of
type $H_4$ in the $4$-dimensional Euclidean space $V$. Let
$\Delta=\{\alpha_1,\alpha_2,\alpha_3,\alpha_4\}$ be a set of
simple roots in $\Phi$ and let
$\{\omega_1,\omega_2,\omega_3,\omega_4\}$ be the corresponding
fundamental weights of $\Delta$. The positive root poset
$(\Phi^+,\leq)$ for type $H_4$ is shown in
Figure~\ref{fi:H4rootposet}, where $\leq$-minimal elements are at
the top and covering relations that do not arise from simple
reflections are marked with dashed lines.  The positive roots for
$H_4$ are given in terms of the basis of simple roots and numbered
consecutively across rows and from top to bottom.

\subsection{Definition}
Let $\Lambda$ be a $\subseteq$-maximal antichain in
$(\Phi^+,\leq)$.
\begin{enumerate}
\item $\Lambda$ is said to be a \textit{good} antichain if
$\mathrm{Int}_C(\Lambda) \neq \emptyset$ and
\item $\Lambda$ is said to be a \textit{bad} antichain if
$\mathrm{Int}_C(\Lambda)=\emptyset$.
\end{enumerate}

\subsection{Remark}

Unlike for types $H_3$ and $I_2(m)$, in type $H_4$ there exist bad
antichains, so the hypothesis in
Corollary~\ref{cor:reduction-thru-antichain} (i) fails to hold for
some $\subseteq$-maximal antichains in type $H_4$.  For example,
for each antichain $I_{\min}$ in Table~\ref{ta:bad4} in the
Appendix, the solution to $(v \mid \beta)=1$ for all $\beta \in
I_{\min}$ does not lie in $C$.  Although each antichain in
Table~\ref{ta:bad4} turns out to define a non-empty region, there
exist other antichains in $(\Phi^+,\leq)$ that do define empty
regions.

Recall that the function $f: \mathcal{R} \rightarrow \mathcal{I}$
given by $f(R)=I_R$ is injective, where
\[
 I_{R}=\{ \beta \in \Phi^+ \mid R \subseteq H_{\beta,1}^+ \}\,.
\]
We will show in general in Section~\ref{se:Criterion} that $f$ is
not surjective if and only if there exist antichains with
$\mathrm{Int}_C(\Lambda)=\emptyset$.  Since $f$ is not surjective
when $\Phi$ is the root system of type $H_4$, we provide instead
an explicit description of the image set of $f$.

\begin{theorem}
\label{thm:H4} For $(\Phi^+,\leq)$ of type $H_4$, the image set
of the function $f: \mathcal{R} \rightarrow \mathcal{I}$, given by
$f: R \mapsto I_R$ for $R \in \mathcal{R}$, is
 $\mathcal{I} \smallsetminus \{\,I(\Lambda_i) \mid 1 \leq i \leq 16\,\}$, where
\[
 \begin{array}{llll}
  \Lambda_1=\{\alpha_{16}\},
 &\Lambda_5=\{\alpha_{14},\alpha_{21}\},
 &\Lambda_9\phantom{i}=\{\alpha_{16},\alpha_{25}\},
 &\Lambda_{13}=\{\alpha_{14},\alpha_{21},\alpha_{23}\}, \\
  \Lambda_2=\{\alpha_{13},\alpha_{16}\},
 &\Lambda_6=\{\alpha_{14},\alpha_{25}\},
 &\Lambda_{10}=\{\alpha_{18},\alpha_{21}\},
 &\Lambda_{14}=\{\alpha_{14},\alpha_{23},\alpha_{25}\}, \\
  \Lambda_3=\{\alpha_{13},\alpha_{20}\},
 &\Lambda_7=\{\alpha_{16},\alpha_{17}\},
 &\Lambda_{11}=\{\alpha_{18},\alpha_{25}\},
 &\Lambda_{15}=\{\alpha_{14},\alpha_{25},\alpha_{28}\}, \\
  \Lambda_4=\{\alpha_{14},\alpha_{19}\},
 &\Lambda_8=\{\alpha_{16},\alpha_{21}\},
 &\Lambda_{12}=\{\alpha_{21},\alpha_{22}\},
 &\Lambda_{16}=\{\alpha_{18},\alpha_{25},\alpha_{28}\}.
 \end{array}
\]
\end{theorem}

\begin{proof}
For type $H_4$, $(\Phi^+,\leq)$ has 60 elements and 429 antichains
(confirmed using the posets package~\cite{St} in Maple), of which
\begin{align*}
20 &\hbox{ are of size four (see Tables~\ref{ta:bad4} and~\ref{ta:good4}
in the Appendix)}, \\
142 &\hbox{ are of size three, and} \\
206 &\hbox{ are of size two}.
\end{align*}

\smallskip
\textbf{Step 1}: Find all $\Lambda$ to which
Corollary~\ref{cor:reduction-thru-antichain} applies.  There are
152 $\subseteq$-maximal antichains. This includes all 20
antichains of size 4, the 79 antichains of size 3 in
Tables~\ref{ta:bad3} and~\ref{ta:good3a}, 47 antichains of size 2,
and the 6 antichains $\{\alpha_i\}$ for $55\leq i \leq 60$. The 12
good antichains of size 4 and the 74 good
antichains of size 3 are listed in Tables~\ref{ta:good4}
and~\ref{ta:good3a}.  By
Proposition~\ref{prop:criterion-for-incomparable}, all 47
$\subseteq$-maximal antichains of size 2 are good and all 6
$\subseteq$-maximal antichains of size 1 are trivially good.

To each of the 139 good antichains, $\Lambda$, apply
Corollary~\ref{cor:reduction-thru-antichain} to eliminate all
$\Lambda'$, where $\Lambda' \preceq \Lambda$, from the list of
antichains in $(\Phi^+,\leq)$.  Note that even if $\Lambda'
\preceq \Lambda$, quite possibly $|\Lambda'| \geq |\Lambda|$.
This process, easily implemented in Maple, yields that 401 of the
antichains define nonempty regions and leaves 28 antichains which
must be investigated by other means.

\smallskip
\textbf{Step 2}: To show $R_{\Lambda_i}=\emptyset$ for $1 \leq i
\leq 16$, suppose instead there exists $v \in V$ such that
$$ (v \mid \beta)>1 \hbox{ for all } \beta \in I_{\min} \hbox{ and }
(v \mid \gamma)<1 \hbox{ for all } \gamma \in I^c_{\max}.$$
Then it is also true that
$\left(v \mid \sum_{j=1}^k c_j \beta_j\right)>1 \hbox{ and }
\left(v \mid \sum_{j=1}^\ell d_j \gamma_j\right)<1$
for any convex linear combinations
\[\sum_{j=1}^k c_j \beta_j \hbox{ and } \sum_{j=1}^\ell d_j \gamma_i
\hbox{ with } \sum_{j=1}^k c_j = \sum_{j=1}^\ell d_j = 1, \beta_j \in I_{\min},
\hbox{ and } \gamma_j \in I^c_{\max}.\] Table~\ref{ta:comparisons}
in the Appendix contains, for each
$I_{\min}=\Lambda_1,\dots,\Lambda_{16}$, the corresponding
$I^c_{\max}$ and a pair of convex linear combinations with $\sum_{j=1}^k c_j \beta_j
< \sum_{j=1}^\ell d_j \gamma_j$.  In each case these provide a contradiction
to the assumption that $R_{\Lambda_i}\neq \emptyset$.  In three cases
the same pair of convex combinations works for more than one $\Lambda$.

\smallskip
\textbf{Step 3}: To verify for the twelve remaining $\Lambda$ that
$R_\Lambda\neq \emptyset$, apply
Theorem~\ref{thm:reduction-thru-antichain}, or a minor
modification of it in the one case of
$\Lambda=\{\alpha_{18},\alpha_{23},\alpha_{25}\}$.  The necessary
calculations are summarized in Tables~\ref{ta:solutionsa}
and~\ref{ta:solutionsb} in the Appendix according to whether part
(i) or (ii) of Theorem~\ref{thm:reduction-thru-antichain} is
used.  For the six antichains which satisfy one of the hypotheses
in Theorem~\ref{thm:reduction-thru-antichain}, the table includes
the solution to the appropriate system of linear equations, and
for the other five antichains, $\Lambda$, the table includes a
note that either $I(\Lambda) = I(\Lambda') \setminus
\Lambda''$ or $I(\Lambda)=I(\Lambda') \cup \Lambda''$, for some
such $\Lambda'$ and $\Lambda'' \subseteq \Lambda'$.

In the one final case of
$I_{\min}=\{\alpha_{18},\alpha_{23},\alpha_{25}\}$ and
 $I_{\max}^c=\{\alpha_{14},\alpha_{21},\alpha_{24}\}$, use the
same notation as in the Appendix with
$(c_1,\,c_2,\,c_3,\,c_4)=c_1\omega_1+c_2\omega_2+c_3\omega_3+c_4\omega_4$.
Then a solution $v \in C$ to $(v \mid \alpha_{18})=(v \mid
\alpha_{25})=1$ has
\[
 v=(c,\,d,\,\tau-1-\tau c-2d,\,(\tau-1)c+d),
\]
with $c>0$, $d>0$, and $\tau c+2d<\tau-1$. Since
$\alpha_{21}<\alpha_{25}$ and $\alpha_{14}<\alpha_{18}$ we know
$(v \mid \alpha_{21})<1$ and $(v \mid \alpha_{14})<1$.
If we take $c=\frac{18-9\tau}{10}$ and $d=\frac{\tau-1}{40}$, then
\begin{align*}
 (v \mid \alpha_{23})&=\tau-1+\tau c+\tau d>1 \hbox{ and} \\
 (v \mid \alpha_{24})&=\tau-2c-\tau d<1.
\end{align*}
Thus by a generalization of
Theorem~\ref{thm:reduction-thru-antichain} which can be proved by
slight modification of the given proof,
$R_{\{\alpha_{18},\alpha_{23},\alpha_{25}\}} \neq \emptyset$.
\end{proof}

\subsection{Remark}
Steps 1 and 3 in the above proof could be combined by applying
Theorem~\ref{thm:reduction-thru-antichain} to all antichains at
the outset, followed by the same arguments as in Step 2 since
$\Lambda_1,\dots,\Lambda_{16}$ would still need to be
investigated.  While more direct, this is computationally more
intensive and we performed the calculations as outlined above.

\subsection{Consequences}
Using the posets for $H_3$ and $H_4$
in Figure~\ref{fi:H4rootposet} and the results in this and the previous
section leads to the numbers of dominant regions shown in
Table~\ref{ta:numbers}.  As described in~\cite[Lemma 4.1]{A2}, a
dominant region, $R$, is bounded if and only if $R$ is contained
in the parallelepiped spanned by $\omega_1, \omega_2, \ldots,
\omega_n$. This in turn occurs if and only if the antichain
corresponding to $R$ does not contain any simple roots.  The proof
of~\cite[Lemma 4.1]{A2} still holds when $\Phi$ is
noncrystallographic, provided that for every
$\beta=\sum\limits_{i=1}^n c_i \alpha_i \in \Phi^+$, all $c_i>1$.
This is true except in type $I_2(m)$ when $m$ is even and $r$
takes on extreme values (see Figure~\ref{fi:I2(6)}). Table 2 lists
the values of the generalized Catalan numbers for comparison.
\begin{table}[h]
\caption{Numbers of dominant regions in noncrystallographic type}
\label{ta:numbers}
\begin{tabular}{|c||c|c|c|c|}
  \hline
  Type & $I_2(m), m \hbox{ odd}$ & $I_2(m), m \hbox{ even}$ & $H_3$ & $H_4$  \\
  \hline \hline
$\begin{array}{ll}
\text{Number of} \\
\text{dominant regions}
\end{array}$ & $\dfrac{3m+1}{2}$ &
$\begin{array}{ll}
    \frac{3m}{2}-1, & \hbox{if $r=\frac{\sin(k\pi/m)}{\sin(\ell\pi/m)}, 1 \leq k,\ell <m/2$;} \\
    \frac{3m}{2}, & \hbox{if $r=\sin(k\pi/m), 1 \leq k <m/2$;} \\
    \frac{3m}{2}+1, & \hbox{otherwise.} \\
\end{array}$
    & 41 & 413 \\
    \hline
$\begin{array}{ll}
\text{Number of bounded} \\
\text{dominant regions}
\end{array}$ & $\dfrac{3m+1}{2}-3$ &
$\begin{array}{ll}
    \frac{3m}{2}-4, & \hbox{if $r=\frac{\sin(k\pi/m)}{\sin(\ell\pi/m)}, 1 \leq k,\ell <m/2$;} \\
    \frac{3m}{2}-3, & \hbox{if $r=\sin(k\pi/m), 1 \leq k <m/2$;} \\
    \frac{3m}{2}-2, & \hbox{otherwise.} \\
\end{array}$
    & 29 & 355 \\
  \hline
\end{tabular}
\end{table}
\begin{table}[!ht]
\caption{Generalized Catalan numbers in noncrystallographic type}
\label{ta:Catalan}
\begin{tabular}{|c||c|c|c|}
  \hline
  Type & $I_2(m)$ & $H_3$ & $H_4$  \\
  \hline \hline
  $e_1,\dots,e_n$ & 1, $m-1$ & 1, 5, 9 & 1, 11, 19, 29 \\  \hline
  $h$ & $m$ & 10 & 30 \\
  \hline
  $\prod\limits_{i=1}^n\frac{e_i+1+h}{e_i+1}$ & $m+2$ & 32 & 280 \\
  \hline
  $\prod\limits_{i=1}^n\frac{e_i-1+h}{e_i-1}$ & $m-1$ & 21 & 232 \\
  \hline
\end{tabular}
\end{table}

\section{Criterion for Empty Regions}
\label{se:Criterion}

In this section we prove the following criterion.
\begin{theorem}  \label{thm:criterion}
The function $f:\mathcal{R} \rightarrow \mathcal{I}$ given by
$f(R)=I_{R}$ for $R \in \mathcal{R}$ is a bijection if and only if
$\Int_{C}(\Lambda) =\{ v \in C \mid (v \mid \beta)=1 \mbox{ for
all } \beta \in \Lambda \} \ne \emptyset$ for any nonempty
antichain $\Lambda$ in $(\Phi^+,\leq)$.
\end{theorem}

The proof of Theorem~\ref{thm:criterion} follows the lemmas
below.  Let $\Lambda$ be a nonempty antichain in $(\Phi^+,\leq)$.
By Proposition~\ref{prop:criterion-for-incomparable}, if
$\Int_{C}(\Lambda)=\emptyset$, then $\Lambda$ must have at least
three positive roots.

\subsection{Notations}
Let $\Lambda$ be an antichain in $(\Phi^+,\leq)$ which has at
least two positive roots.
\begin{enumerate}
\item  Given any $\gamma \in \Lambda$, define
$\Lambda_{\gamma}=\Lambda \smallsetminus \{ \gamma \}$.

\item  Let $\mathcal{A}(\Lambda)$ be the set of all the connected
components of
 $V \smallsetminus \bigcup_{\beta \in \Lambda} H_{\beta,1}$
under the Euclidean topology.
 Any region $A \in \mathcal{A}(\Lambda)$
can be uniquely expressed as
 $\bigcap_{\beta \in \Lambda} H_{\beta,1}^{n_{\beta}^A}$,
where $n_{\beta}^A$ represents $+$ or $-$. So we can write $A$
simply as $(n_{\beta}^A)_{\beta \in \Lambda}$ (see \cite{Sh2}).

\item  For any distinct vectors $u,\,v \in V$, we define
$l_{u,v}(t)=(1-t)u+tv$ for $t \in \mathbb{R}$.
\end{enumerate}

Using Proposition~\ref{prop:linear-independence-of-antichain} it
is easy to show that given a $\Lambda$-tuple $(n_{\beta})_{\beta
\in \Lambda}$ with $n_{\beta}=+$ or $-$,
\begin{equation}
 \tag{**} \label{lem:tuple-and-region}
 \hbox{there exists a unique region $A \in \mathcal{A}(\Lambda)$
 such that $n_{\beta}^A=n_{\beta}$ for all $\beta \in \Lambda$.}
\end{equation}

\begin{lemma}  \label{lem:intersection-with-convex-hull}
Let $\gamma \in \Lambda$ and let $v_A \in A \cap H_{\gamma,1}$ for
each $A \in \mathcal{A}(\Lambda_{\gamma})$. Then
\[ \Int(\Lambda) \cap
   \Conv\{ v_A \}_{A \in \mathcal{A}(\Lambda_{\gamma})}
   \ne \emptyset\,, \]
where
 $\Conv\{ v_A \}_{A \in \mathcal{A}(\Lambda_{\gamma})}$
is the convex hull in $V$ generated by $v_A$ for $A \in
\mathcal{A}(\Lambda_{\gamma})$.
\end{lemma}

\begin{proof}
If $|\Lambda|=2$, then $\Lambda_{\gamma}=\{\beta\}$ and
$\mathcal{A}(\Lambda_{\gamma})=\{A^+,A^-\}$, where
$A^+=H_{\beta,1}^+$ and $A^-=H_{\beta,1}^-$. Set
 $t^*=\frac{(v_{A^+} \mid \beta)-1}
      {(v_{A^+} \mid \beta)-(v_{A^-} \mid \beta)}$.
Then $0<t^*<1$ and
 $l_{v_{A^+},v_{A^-}}(t^*) \in \Int(\Lambda) \cap
 \Conv\{ v_{A^+}, v_{A^-}\}$.

Suppose that this lemma holds for all sets $\Lambda'$ with
$\{\gamma\} \subsetneq \Lambda' \subsetneq \Lambda$. Take
 $\beta_1 \in \Lambda_{\gamma}$ and set
 $\Lambda'=\Lambda \smallsetminus \{\beta_1\}$.
Then every region $A \in \mathcal{A}(\Lambda_{\gamma})$ can be
expressed as a $\Lambda_{\gamma}$-tuple
 $(n_{\beta_1}^A, (n_{\beta}^A)_{\beta \in \Lambda_{\gamma}'})$
and by~(\ref{lem:tuple-and-region}) there is some region
 $A' \in \mathcal{A}(\Lambda_{\gamma}')$
such that
 $n_{\beta}^{A'}=n^{A}_{\beta}$
for any $\beta \in \Lambda_{\gamma}'$.
 Every  region $A' \in \mathcal{A}(\Lambda_{\gamma}')$
is partitioned into two disjoint subregions $A_{A'}^+$ and
$A_{A'}^-$ in $\mathcal{A}(\Lambda_{\gamma})$ which are
represented by the $\Lambda_{\gamma}$-tuples
 $(+,(n^{A'}_{\beta})_{\beta \in \Lambda_{\gamma}'})$ and
 $(-,(n^{A'}_{\beta})_{\beta \in \Lambda_{\gamma}'})$,
respectively. Since $|\Lambda'|<|\Lambda|$, we can find
 $v^+ \in \Int(\Lambda') \cap
   \Conv\{ v_{A_{A'}^+} \}_{A' \in \mathcal{A}(\Lambda_{\gamma}')}$
and
 $v^- \in \Int(\Lambda') \cap
   \Conv\{ v_{A_{A'}^-} \}_{A' \in \mathcal{A}(\Lambda_{\gamma}')}$.
Set
 $t^*=\frac{(v^+ \mid \beta_1)-1}
      {(v^+ \mid \beta_1)-(v^- \mid \beta_1)}$.
Then $0<t^*<1$ and
 $l_{v^+,v^-}(t^*) \in \Int(\Lambda) \cap
  \Conv\{ v_A \}_{A \in \mathcal{A}(\Lambda_\gamma)}$.
So the lemma follows by induction on $|\Lambda|$.
\end{proof}

\begin{lemma}  \label{lem:empty-intersection-of-regions}
If $\Int_{C}(\Lambda)=\emptyset$ and
   $\Int_{C}(\Lambda_{\gamma}) \ne \emptyset$
for some $\gamma \in \Lambda$, then there exists a nonempty region
 $B=\bigcap_{\beta \in \Lambda_{\gamma}}
      (H_{\beta,1}^{n_\beta} \cap C)$
such that
 $B \cap H_{\gamma,1}=\emptyset$,
where $n_{\beta}=+$ or $-$ for each $\beta \in \Lambda_{\gamma}$.
\end{lemma}

\begin{proof}
Since $\Int_{C}(\Lambda_{\gamma}) \ne \emptyset$ and
$\Lambda_{\gamma}$ is a linearly independent set in $\Lambda$,
$\bigcap_{\beta \in \Lambda_{\gamma}}
  (H_{\beta,1}^{n_{\beta}} \cap C)
  \ne \emptyset$
for any choice of $n_{\beta}$ by an argument similar to that in
step 3 of the proof of Theorem~\ref{thm:reduction-thru-antichain}.

Suppose that the desired region in this lemma does not exist. For
each region $A \in \mathcal{A}(\Lambda_{\gamma})$, we choose
 $v_A \in (A \cap C) \cap H_{\gamma,1}$. By
Lemma~\ref{lem:intersection-with-convex-hull}, we obtain
 $\Int(\Lambda) \cap
  \Conv\{v_{A}\}_{A \in \mathcal{A}(\Lambda_{\gamma})}
  \ne \emptyset$,
which yields
 $\Int(\Lambda_{\gamma}) \cap
  \Conv\{v_{A}\}_{A \in \mathcal{A}(\Lambda_{\gamma})}
  \ne \emptyset$.
Since
 $\Conv\{ v_{A} \}_{A \in \mathcal{A}(\Lambda_{\gamma})}
  \subseteq C \cap H_{\gamma,1}$,
we get
\[ \Int_{C}(\Lambda)
  =\Int_{C}(\Lambda_{\gamma}) \cap H_{\gamma,1}
  \supseteq
   \Int(\Lambda_{\gamma}) \cap
   \Conv\{ v_{A} \}_{A \in \mathcal{A}(\Lambda_{\gamma})}
  \ne \emptyset\,, \]
contrary to the hypothesis $\Int_{C}(\Lambda)=\emptyset$.
\end{proof}

\begin{lemma}  \label{lem:criterion-for-empty-regions}
If $\Int_{C}(\Lambda)=\emptyset$ for some nonempty antichain
$\Lambda$ in $(\Phi^+,\leq)$, then there exists an increasing set
 $I \in \mathcal{I}$ such that $R_I=\emptyset$.
\end{lemma}

\begin{proof}
Without loss of generality, assume that there is some $\gamma \in
\Lambda$ such that
  $\Int_C(\Lambda_{\gamma}) \ne \emptyset$.
By Lemma \ref{lem:empty-intersection-of-regions}, there exists a
nonempty region
 $B=\bigcap_{\beta \in \Lambda_{\gamma}}
 (H_{\beta,1}^{n_{\beta}} \cap C)$
such that $B \cap H_{\gamma,1}=\emptyset$, where $n_{\beta}=+$ or
$-$ for each $\beta \in \Lambda_{\gamma}$. Let
 $\Lambda_1=\{\beta \in \Lambda_{\gamma} \mid n_{\beta}=+\}$,
which is an antichain in $(\Phi^+,\leq)$, and set
$I_1=I(\Lambda_1)$, which is the increasing set in $\mathcal{I}$
generated by $\Lambda_1$. If $R_{I_1}=\emptyset$, then we are
done. Now, we assume that $R_{I_1} \ne \emptyset$. Since $\Lambda$
is an antichain, then
  $\Lambda \smallsetminus \Lambda_1
  \subseteq \Phi^+ \smallsetminus I_1$.
In particular, given any $v \in R_{I_1}$ we have
 $(v \mid \beta)<1$
for all $\beta \in \Lambda_{\gamma} \smallsetminus \Lambda_1$.
Thus, $R_{I_1} \subseteq B$. Note that $B$ is a (connected) convex
region in $V$ and $B \cap H_{\gamma,1}=\emptyset$. There are
exactly two cases to be investigated:

\textit{Case 1:} If $B \subseteq H_{\gamma,1}^+$, then
 $(v \mid \gamma)>1$ for $v \in R_{I_1}$, contrary to the fact that
  $\gamma \not\in I_1$.

\textit{Case 2:} If $B \subseteq H_{\gamma,1}^-$, then set
$\Lambda_2=\Lambda_1 \cup \{\gamma\}$, which is an antichain in
$(\Phi^+,\leq)$, and $I_2=I(\Lambda_2)$, which is the increasing
set in $\mathcal{I}$ generated by $\Lambda_2$. Suppose that
$R_{I_2} \ne \emptyset$. Then we still have
 $R_{I_2} \subseteq B \subseteq H_{\gamma,1}^-$.
By $\gamma \in I_2$, we would have $(v \mid \gamma)>1$ for all $v
\in R_{I_2}$, contrary to the fact that $R_{I_2} \subseteq
H_{\gamma,1}^-$.
\end{proof}

\subsection{Proof of Theorem~\ref{thm:criterion}}
\quad

$\Rightarrow)$  Suppose that $\Int_{C}(\Lambda)=\emptyset$ for
some nonempty antichain $\Lambda$ in $(\Phi^+,\leq)$. By
Lemma~\ref{lem:criterion-for-empty-regions}, we can find an
increasing set $I \in \mathcal{I}$ such that $R_I=\emptyset$.
Therefore, $I$ does not lie in the image set of $f$ by
Lemma~\ref{lem:region-to-increasing set}.

$\Leftarrow)$  Suppose that $\Int_{C}(\Lambda) \ne \emptyset$ for
any nonempty antichain $\Lambda$ in $(\Phi^+,\leq)$. Let $I$ be a
nonempty increasing set in $(\Phi^+,\leq)$. Using the same
argument as in the proof of Theorem~\ref{thm:main}, we can show
$R_I \ne \emptyset$. Therefore we get $R_I \in \mathcal{R}$ which
satisfies $f(R_I)=I_{R_I}=I$.  \qed

\vfill \newpage
\bigskip
\begin{appendix}

\section*{Appendix. Calculations used in the proof of Theorem~\ref{thm:H4}}
For $\Lambda$ a nonempty antichain
in $(\Phi^+,\leq)$, the notation $(c_1,c_2,c_3,c_4)$ in the tables
in this section represents a solution
$c_{1}\omega_{1}+c_{2}\omega_{2}+c_{3}\omega_{3}+c_{4}\omega_{4}$
to the system $(v \mid \beta)=1$ with $\beta \in \Lambda$.
\vspace{-.1in}
\begin{table}[!ht]
\caption{Bad $4$-element antichains in type $H_4$}
\label{ta:bad4}
\begin{tabular}{|c|c|}
\hline
  $I_{\min}$ & Solution to $(v \mid \beta)=1$ for all $\beta \in I_{\min}$
\\ \hline\hline
  $\{\alpha_{13}, \alpha_{14}, \alpha_{16}, \alpha_{19}\}$
 & $(2-\tau,\, 3\tau-5,\, 2-\tau,\, 5-3\tau)$
\\ \hline

  $\{\alpha_{13}, \alpha_{14}, \alpha_{19}, \alpha_{20}\}$
 & $(0,\, 2\tau-3,\, 2-\tau,\, 0)$
\\ \hline

  $\{\alpha_{14}, \alpha_{16}, \alpha_{17}, \alpha_{19}\}$
 & $(2-\tau,\, 3-2\tau,\, \tau-1,\, 0)$
\\ \hline

  $\{\alpha_{14}, \alpha_{21}, \alpha_{23}, \alpha_{24}\}$
 & $(1,\, -1,\, \tau-1,\, 2-\tau)$
\\ \hline

  $\{\alpha_{14}, \alpha_{23}, \alpha_{24}, \alpha_{25}\}$
 & $(1-\tau,\, 1,\, \tau-2,\, 2-\tau)$
\\ \hline

  $\{\alpha_{18}, \alpha_{21}, \alpha_{23}, \alpha_{24}\}$
 & $(2-\tau,\, 3\tau-5,\, 2\tau-3,\, 5-3\tau)$
\\ \hline

  $\{\alpha_{18}, \alpha_{23}, \alpha_{24}, \alpha_{25}\}$
 & $(\tau-1,\, \tau-2,\, 2-\tau,\, 0)$
\\ \hline

  $\{\alpha_{21}, \alpha_{22}, \alpha_{23}, \alpha_{24}\}$
 & $(0,\, 2-\tau,\, 0,\, 0)$
\\ \hline
\end{tabular}
\end{table}
\vspace{-.1in}
\begin{table}[!ht]
\caption{Good $4$-element antichains in type $H_4$}
\label{ta:good4}
\begin{tabular}{|c|c|}
\hline
  $I_{\min}$ & Solution to $(v \mid \beta)=1$ for all $\beta \in I_{\min}$
\\ \hline\hline
  $\{\alpha_1, \alpha_2, \alpha_3, \alpha_4\}$
 & $(1,\, 1,\, 1,\, 1)$
\\ \hline

  $\{\alpha_3, \alpha_4, \alpha_5, \alpha_6\}$
 & $(2-\tau,\, 2-\tau,\, 1,\, 1)$
\\ \hline

  $\{\alpha_4, \alpha_5, \alpha_6, \alpha_7\}$
 & $(2-\tau,\, 2-\tau,\, \tau-1,\, 1)$
\\ \hline

  $\{\alpha_4, \alpha_9, \alpha_{10}, \alpha_{11}\}$
 & $(2-\tau,\, 2\tau-3,\, 5-3\tau,\, 1)$
\\ \hline

  $\{\alpha_5, \alpha_6, \alpha_7, \alpha_8\}$
 & $(2-\tau,\, 2-\tau,\, \tau-1,\, 2-\tau)$
\\ \hline

  $\{\alpha_8, \alpha_9, \alpha_{10}, \alpha_{11}\}$
 & $(2-\tau,\, 2\tau-3,\, 5-3\tau,\, 3\tau-4)$
\\ \hline

  $\{\alpha_9, \alpha_{10}, \alpha_{11}, \alpha_{12}\}$
 & $(2-\tau,\, 2\tau-3,\, 5-3\tau,\, \tau-1)$
\\ \hline

  $\{\alpha_{13}, \alpha_{14}, \alpha_{15}, \alpha_{16}\}$
 & $(2-\tau,\, 5\tau-8,\, 5-3\tau,\, 5-3\tau)$
\\ \hline

  $\{\alpha_{14}, \alpha_{17}, \alpha_{19}, \alpha_{20}\}$
 & $(5-3\tau,\, 5-3\tau,\, 2\tau-3,\, 5-3\tau)$
\\ \hline

  $\{\alpha_{17}, \alpha_{18}, \alpha_{19}, \alpha_{20}\}$
 & $\frac{1}{2}(2\tau-3,\, 2-\tau,\, 2-\tau,\, 2\tau-3)$
\\ \hline

  $\{\alpha_{22}, \alpha_{23}, \alpha_{24}, \alpha_{25}\}$
 & $(5-3\tau,\, 5-3\tau,\, 5\tau-8,\, 5-3\tau)$
\\ \hline

  $\{\alpha_{25}, \alpha_{26}, \alpha_{27}, \alpha_{28}\}$
 & $(13-8\tau,\, 2\tau-3,\, 13-8\tau,\, 2-\tau)$
\\ \hline
\end{tabular}
\end{table}
\vspace{-.1in}
\begin{table}[!ht]
\caption{Bad $3$-element antichains in type $H_4$}
\label{ta:bad3}
\begin{tabular}{|c|r@{\,+\,$c$\,}l|}
\hline
  $I_{\min}$ &
\multicolumn{2}{c|}{
  Solutions to $(v \mid \beta)=1$ for all $\beta \in I_{\min}$}
\\ \hline\hline
  $\{\alpha_{14},\alpha_{16},\alpha_{21}\}$
 & $(2-\tau,\,0,\,0,\,2-\tau)$&$(0,\,1,\,\tau-2,\,\tau-1)$
\\ \hline

  $\{\alpha_{14},\alpha_{16},\alpha_{25}\}$
 & $(2-\tau,\,0,\,0,\,2-\tau)$&$(0,\,1,\,-2,\,1)$
\\ \hline

  $\{\alpha_{14},\alpha_{19},\alpha_{24}\}$
 & $(2-\tau,\,0,\,0,\,2-\tau)$&$(\tau,\,-2,\,1,\,0)$
\\ \hline

  $\{\alpha_{14},\alpha_{25},\alpha_{28}\}$
 & $\frac{1}{5}(4-3\tau,\,2\tau-1,\,0,\,3-\tau)$&$
 (-3\tau-1,\,2\tau-1,\,1,\,-\tau-2)$
\\ \hline

  $\{\alpha_{18},\alpha_{25},\alpha_{28}\}$
 & $(\tau-2,\,\tau-1,\,0,\,2-\tau)$&$(-\tau-1,\,\tau,\,1,\,0)$
\\ \hline
\end{tabular}
\end{table}

\begin{table}[!ht]
\caption{Good $3$-element antichains in type $H_4$}
\label{ta:good3a}
\begin{tabular}{|c|r@{\,+\,$c$\,}l|}
\hline
  $I_{\min}$ &
\multicolumn{2}{c|}{
  Solutions to $(v \mid \beta)=1$ for all $\beta \in I_{\min}$}
\\ \hline\hline
   $\{\alpha_1,\alpha_2,\alpha_8\}$
 & $\frac{1}{2}(2,\,2,\,1,\,1)$&$(0,\,0,\,1,\,-1)$
\\ \hline

   $\{\alpha_1,\alpha_4,\alpha_7\}$
 & $\frac{1}{2}(2,\,1,\,1,\,2)$&$(0,\,1,\,-1,\,0)$
\\ \hline

   $\{\alpha_1,\alpha_7,\alpha_8\}$
 & $\frac{1}{2}(2,\,1,\,1,\,1)$&$(0,\,1,\,-1,\,1)$
\\ \hline

   $\{\alpha_3,\alpha_4,\alpha_9\}$
 & $\frac{1}{2}(\tau-1,\,\tau-1,\,2,\,2)$&$(1,\,-1,\,0,\,0)$
\\ \hline

   $\{\alpha_4,\alpha_5,\alpha_{11}\}$
 & $\frac{1}{2}(2-\tau,\,1,\,2-\tau,\,2)$&$(-1,\,\tau-1,\,1,\,0)$
\\ \hline

   $\{\alpha_4,\alpha_6,\alpha_{10}\}$
 & $\frac{1}{2}(1,\,2-\tau,\,2\tau-3,\,2)$&$(1,\,-\tau,\,1,\,0)$
\\ \hline

   $\{\alpha_4,\alpha_7,\alpha_9\}$
 & $\frac{1}{2}(\tau-1,\,\tau-1,\,3-\tau,\,2)$&$(1,\,-1,\,1,\,0)$
\\ \hline

   $\{\alpha_4,\alpha_{10},\alpha_{15}\}$
 & $\frac{1}{2}(2-\tau,\,2\tau-3,\,4-2\tau,\,2)$&$(\tau,\,-1,\,0,\,0)$
\\ \hline

   $\{\alpha_4,\alpha_{10},\alpha_{19}\}$
 & $\frac{1}{4}(2-\tau,\,4\tau-6,\,5-2\tau,\,4)$&$(\tau,\,-2,\,1,\,0)$
\\ \hline

   $\{\alpha_4,\alpha_{13},\alpha_{15}\}$
 & $\frac{1}{2}(\tau-1,\,2\tau-3,\,2-\tau,\,2)$&$(-\tau,\,\tau-1,\,1,\,0)$
\\ \hline

   $\{\alpha_4,\alpha_{13},\alpha_{19}\}$
 & $\frac{1}{2}(2\tau-3,\,2\tau-3,\,4-2\tau,\,2)$&$(1\,-1,\,0,\,0)$
\\ \hline

   $\{\alpha_4,\alpha_{17},\alpha_{19}\}$
 & $\frac{1}{2}(2\tau-3,\,2-\tau,\,2-\tau,\,2)$&$(\tau-1,\,-1,\,1,\,0)$
\\ \hline

   $\{\alpha_5,\alpha_6,\alpha_{12}\}$
 & $\frac{1}{2}(4-2\tau,\,4-2\tau,\,\tau-1,\,\tau-1)$&$(0,\,0,\,1,\,-1)$
\\ \hline

   $\{\alpha_5,\alpha_8,\alpha_{11}\}$
 & $\frac{1}{2}(2-\tau,\,1,\,2-\tau,\,\tau)$&$(-1,\,\tau-1,\,1,\,-1)$
\\ \hline

   $\{\alpha_5,\alpha_{11},\alpha_{12}\}$
 & $\frac{1}{2}(2-\tau,\,1,\,2-\tau,\,\tau-1)$&$(1,\,1-\tau,\,-1,\,\tau)$
\\ \hline

   $\{\alpha_6,\alpha_8,\alpha_{10}\}$
 & $\frac{1}{2}(1,\,2-\tau,\,2\tau-3,\,5-2\tau)$&$(1,\,-\tau,\,1,\,-1)$
\\ \hline

   $\{\alpha_6,\alpha_{10},\alpha_{12}\}$
 & $\frac{1}{2}(1,\,2-\tau,\,2\tau-3,\,3-\tau)$&$(1,\,-\tau,\,1,\,\tau-1)$
\\ \hline

   $\{\alpha_7,\alpha_8,\alpha_9\}$
 & $\frac{1}{2}(\tau-1,\,\tau-1,\,3-\tau,\,\tau-1)$&$(1,\,-1,\,1,\,-1)$
\\ \hline

   $\{\alpha_8,\alpha_{10},\alpha_{15}\}$
 & $\frac{1}{2}(2-\tau,\,2\tau-3,\,4-2\tau,\,2\tau-2)$&$(-\tau,\,1,\,0,\,0)$
\\ \hline

   $\{\alpha_8,\alpha_{10},\alpha_{19}\}$
 & $\frac{1}{4}(2-\tau,\,4\tau-6,\,5-2\tau,\,2\tau-1)$&$(\tau,\,-2,\,1,\,-1)$
\\ \hline

   $\{\alpha_8,\alpha_{13},\alpha_{15}\}$
 & $\frac{1}{2}(\tau-1,\,2\tau-3,\,2-\tau,\,\tau)$&$(-\tau,\,\tau-1,\,1,\,-1)$
\\ \hline

   $\{\alpha_8,\alpha_{13},\alpha_{19}\}$
 & $\frac{1}{2}(2\tau-3,\,2\tau-3,\,4-2\tau,\,2\tau-2)$&$(1,\,-1,\,0,\,0)$
\\ \hline

   $\{\alpha_8,\alpha_{17},\alpha_{19}\}$
 & $\frac{1}{2}(2\tau-3,\,2-\tau,\,2-\tau,\,\tau)$&$(\tau-1,\,-1,\,1,\,-1)$
\\ \hline

   $\{\alpha_9,\alpha_{10},\alpha_{16}\}$
 & $\frac{1}{2}(2-\tau,\,3\tau-4,\,5-3\tau,\,2-\tau)$&$(1,\,-1,\,2-\tau,\,-1)$
\\ \hline

   $\{\alpha_9,\alpha_{11},\alpha_{14}\}$
 & $\frac{1}{2}(1,\,2\tau-3,\,5-3\tau,\,2\tau-3)$&$(-\tau,\,\tau,\,1,\,-\tau)$
\\ \hline

   $\{\alpha_9,\alpha_{14},\alpha_{16}\}$
 & $\frac{1}{2}(4-2\tau,\,4\tau-6,\,5-3\tau,\,5-3\tau)$&$(0,\,0,\,1,\,-1)$
\\ \hline

   $\{\alpha_{10},\alpha_{12},\alpha_{15}\}$
 & $\frac{1}{2}(2-\tau,\,2\tau-3,\,4-2\tau,\,1)$&$(-\tau,\,1,\,0,\,-1)$
\\ \hline

   $\{\alpha_{10},\alpha_{12},\alpha_{19}\}$
 & $\frac{1}{4}(2-\tau,\,4\tau-6,\,5-2\tau,\,5-2\tau)$&$(\tau,\,-2,\,1,\,1)$
\\ \hline

   $\{\alpha_{10},\alpha_{15},\alpha_{16}\}$
 & $\frac{1}{2}(2-\tau,\,2\tau-3,\,4-2\tau,\,2-\tau)$&$(\tau,\,-1,\,0,\,-\tau)$
\\ \hline

   $\{\alpha_{10},\alpha_{16},\alpha_{19}\}$
 & $\frac{1}{4}(2-\tau,\,4\tau-6,\,5-2\tau,\,6-3\tau)$&$(\tau,\,-2,\,1,\,-\tau)$
\\ \hline

   $\{\alpha_{10},\alpha_{19},\alpha_{20}\}$
 & $\frac{1}{4}(2-\tau,\,4\tau-6,\,5-2\tau,\,2-\tau)$&$(\tau,\,-2,\,1,\,\tau)$
\\ \hline

   $\{\alpha_{12},\alpha_{13},\alpha_{15}\}$
 & $\frac{1}{2}(\tau-1,\,2\tau-3,\,2-\tau,\,3-\tau)$&$(\tau,\,1-\tau,\,-1,\,\tau)$
\\ \hline

   $\{\alpha_{12},\alpha_{13},\alpha_{19}\}$
 & $\frac{1}{2}(2\tau-3,\,2\tau-3,\,4-2\tau,\,1)$&$(1,\,-1,\,0,\,1)$
\\ \hline

   $\{\alpha_{12},\alpha_{17},\alpha_{19}\}$
 & $\frac{1}{2}(2\tau-3,\,2-\tau,\,2-\tau,\,2\tau-2)$&$(\tau-1,\,-1,\,1,\,0)$
\\ \hline

   $\{\alpha_{13},\alpha_{14},\alpha_{23}\}$
 & $\frac{1}{4}(2-\tau,\,4\tau-6,\,\tau,\,5-3\tau)$&$(1,\,2-2\tau,\,2\tau-3,\,2-\tau)$
\\ \hline
\end{tabular}
\end{table}
\setcounter{table}{5}
\begin{table}[!ht]
\caption{Good $3$-element antichains in type $H_4$, continued}
\begin{tabular}{|c|r@{\,+\,$c$\,}l|}
\hline
  $I_{\min}$ &
\multicolumn{2}{c|}{
  Solutions to $(v \mid \beta)=1$ for all $\beta \in I_{\min}$}
\\ \hline\hline

   $\{\alpha_{14},\alpha_{17},\alpha_{23}\}$
 & $\frac{1}{4}(2-\tau,\,4\tau-6,\,7\tau-10,\,15-9\tau)$&$(1,\,2-2\tau,\,1,\,\tau-2)$
\\ \hline

   $\{\alpha_{14},\alpha_{20},\alpha_{21}\}$
 & $\frac{1}{2}(2-\tau,\, 2\tau-3,\, 2\tau-3,\, 7-4\tau)$&$(1,\,1-\tau,\,1-\tau,\,\tau-1)$
\\ \hline

   $\{\alpha_{14},\alpha_{20},\alpha_{25}\}$
 & $\frac{1}{2}(2-\tau,\,2\tau-3,\,5-3\tau,\,\tau-1)$&$(-\tau,\,1,\,\tau-1,\,1-\tau)$
\\ \hline

   $\{\alpha_{17},\alpha_{18},\alpha_{23}\}$
 & $\frac{1}{2}(2\tau-3,\,2\tau-3,\,\tau-1,\,5-3\tau)$&$(1,\,-1,\,\tau-1,\,1-\tau)$
\\ \hline

   $\{\alpha_{18},\alpha_{19},\alpha_{24}\}$
 & $\frac{1}{10}(16\tau-23,\,25-15\tau,\,7-4\tau,\,13\tau-19)$&$(1,\,\tau-3,\,1,\,1-\tau)$
\\ \hline

   $\{\alpha_{18},\alpha_{20},\alpha_{21}\}$
 & $\frac{1}{2}(2\tau-3,\,2-\tau,\,2\tau-3,\,2-\tau)$&$(-1,\,2-\tau,\,1,\,\tau-2)$
\\ \hline

   $\{\alpha_{18},\alpha_{20},\alpha_{25}\}$
 & $\frac{1}{10}(3\tau-4,\,13\tau-19,\,25-15\tau,\,9\tau-12)$&$(-\tau-1,\,1,\,2\tau-1,\,1-\tau)$
\\ \hline

   $\{\alpha_{19},\alpha_{22},\alpha_{24}\}$
 & $\frac{1}{10}(6\tau-8,\,10-5\tau,\,7-4\tau,\,3\tau-4)$&$(2\tau,\,-\tau-2,\,1,\,\tau)$
\\ \hline

   $\{\alpha_{21},\alpha_{23},\alpha_{27}\}$
 & $\frac{1}{2}(10-6\tau,\,2-\tau,\,5-3\tau,\,5\tau-8)$&$(2,\,-\tau-1,\,1,\,\tau-1)$
\\ \hline

   $\{\alpha_{22},\alpha_{25},\alpha_{28}\}$
 & $\frac{1}{10}(25\tau-40,\,9\tau-12,\,7\tau-11,\,32-19\tau)$&$(1-2\tau,\,1,\,\tau,\,-1)$
\\ \hline

   $\{\alpha_{23},\alpha_{25},\alpha_{27}\}$
 & $\frac{1}{2}(2-\tau,\,2\tau-3,\,5-3\tau,\,2\tau-3)$&$(3-\tau,\,-\tau,\,1,\,\tau-2)$
\\ \hline

   $\{\alpha_{24},\alpha_{25},\alpha_{26}\}$
 & $\frac{1}{10}(10\tau-15,\,3-\tau,\,17\tau-26,\,7-4\tau)$&$(4-3\tau,\,1,\,1-\tau,\,2-\tau)$
\\ \hline

   $\{\alpha_{26},\alpha_{27},\alpha_{31}\}$
 & $\frac{1}{2}(5-3\tau,\,5-3\tau,\,5-3\tau,\,\tau-1)$&$(1,\,-1,\,1,\,-1)$
\\ \hline

   $\{\alpha_{26},\alpha_{27},\alpha_{34}\}$
 & $\frac{1}{2}(5\tau-8,\,5-3\tau,\,5\tau-8,\,4-2\tau)$&$(1,\,-\tau,\,1,\,0)$
\\ \hline

   $\{\alpha_{26},\alpha_{27},\alpha_{37}\}$
 & $\frac{1}{10}(25\tau-40,\,7-4\tau,\,25\tau-40,\,9-3\tau)$&$(4-3\tau,\,1,\,4-3\tau,\,\tau-1)$
\\ \hline

   $\{\alpha_{26},\alpha_{27},\alpha_{40}\}$
 & $\frac{1}{20}(25-15\tau,\,14-8\tau,\,25-15\tau,\,19\tau-22)$&$(5,\,-4-2\tau,\,5,\,\tau-3)$
\\ \hline

   $\{\alpha_{27},\alpha_{28},\alpha_{29}\}$
 & $\frac{1}{2}(5-3\tau,\,2\tau-3,\,10-6\tau,\,2\tau-3)$&$(1,\,-\tau,\,0,\,\tau)$
\\ \hline

   $\{\alpha_{27},\alpha_{29},\alpha_{31}\}$
 & $\frac{1}{4}(10-6\tau,\,2\tau-3,\,15-9\tau,\,5\tau-7)$&$(2,\,-\tau,\,1,\,-1)$
\\ \hline

   $\{\alpha_{27},\alpha_{29},\alpha_{34}\}$
 & $\frac{1}{4}(10\tau-16,\,2\tau-3,\,7\tau-11,\,3-\tau)$&$(2\tau+2,\,-3\tau-2,\,1,\,2\tau+1)$
\\ \hline

   $\{\alpha_{27},\alpha_{29},\alpha_{37}\}$
 & $\frac{1}{2}(5\tau-8,\,5\tau-8,\,5-3\tau,\,4-2\tau)$&$(\tau+1,\,-\tau-1,\,1,\,0)$
\\ \hline

   $\{\alpha_{27},\alpha_{29},\alpha_{40}\}$
 & $\frac{1}{12}(15-9\tau,\,4\tau-6,\,25-15\tau,\,13\tau-16)$&$(3,\,-2\tau,\,1,\,\tau-1)$
\\ \hline

   $\{\alpha_{29},\alpha_{30},\alpha_{31}\}$
 & $\frac{1}{4}(4-2\tau,\,2\tau-3,\,5-3\tau,\,\tau-1)$&$(-2\tau-2,\,\tau,\,1,\,2\tau+1)$
\\ \hline

   $\{\alpha_{29},\alpha_{30},\alpha_{34}\}$
 & $\frac{1}{4}(10-6\tau,\,2\tau-3,\,5-3\tau,\,3-\tau)$&$(-2,\,\tau,\,1,\,-1)$
\\ \hline

   $\{\alpha_{29},\alpha_{30},\alpha_{37}\}$
 & $\frac{1}{2}(5-3\tau,\,5\tau-8,\,13-8\tau,\,4-2\tau)$&$(-\tau-1,\,\tau,\,1,\,0)$
\\ \hline

   $\{\alpha_{29},\alpha_{30},\alpha_{40}\}$
 & $\frac{1}{6}(15-9\tau,\,2\tau-3,\,5-3\tau,\,4-\tau)$&$(-3,\,\tau,\,1,\,\tau-1)$
\\ \hline

   $\{\alpha_{32},\alpha_{33},\alpha_{34}\}$
 & $\frac{1}{4}(10-6\tau,\,2\tau-3,\,2-\tau,\,4-2\tau)$&$(2\tau-4,\,\tau-1,\,1,\,-2)$
\\ \hline

   $\{\alpha_{32},\alpha_{33},\alpha_{37}\}$
 & $\frac{1}{2}(13-8\tau,\,5-3\tau,\,2\tau-3,\,2-\tau)$&$(-1,\,\tau-1,\,1,\,-\tau)$
\\ \hline

   $\{\alpha_{32},\alpha_{33},\alpha_{40}\}$
 & $\frac{1}{10}(25-15\tau,\,7\tau-11,\,7-4\tau,\,13-6\tau)$&$(\tau-3,\,\tau-1,\,1,\,-1)$
\\ \hline

   $\{\alpha_{32},\alpha_{36},\alpha_{37}\}$
 & $\frac{1}{2}(5\tau-8,\,5\tau-8,\,2-\tau,\,5-3\tau)$&$(1,\,-1,\,1,\,-\tau)$
\\ \hline

   $\{\alpha_{32},\alpha_{36},\alpha_{40}\}$
 & $\frac{1}{10}(65-40\tau,\,7\tau-11,\,32-19\tau,\,9\tau-12)$&$(3-\tau,\,-\tau,\,1,\,-1)$
\\ \hline

   $\{\alpha_{32},\alpha_{39},\alpha_{40}\}$
 & $\frac{1}{4}(26-16\tau,\,5\tau-8,\,20-12\tau,\,8\tau-12)$&$(2,\,-\tau,\,0,\,0)$
\\ \hline

   $\{\alpha_{35},\alpha_{36},\alpha_{37}\}$
 & $\frac{1}{2}(5\tau-8,\,5-3\tau,\,2\tau-3,\,2\tau-3)$&$(2-\tau,\,1-\tau,\,1,\,-1)$
\\ \hline

   $\{\alpha_{35},\alpha_{36},\alpha_{40}\}$
 & $\frac{1}{2}(13-8\tau,\,5-3\tau,\,5-3\tau,\,4\tau-6)$&$(1,\,-\tau-1,\,\tau+1,\,0)$
\\ \hline

   $\{\alpha_{35},\alpha_{39},\alpha_{40}\}$
 & $\frac{1}{4}(26-16\tau,\,5-3\tau,\,7-4\tau,\,3\tau-4)$&$(2,\,-\tau-1,\,1,\,\tau)$
\\ \hline

   $\{\alpha_{38},\alpha_{39},\alpha_{40}\}$
 & $\frac{1}{4}(4\tau-6,\,5-3\tau,\,2\tau-3,\,2-\tau)$&$(-2,\,\tau-1,\,1,\,\tau)$
\\ \hline
\end{tabular}
\end{table}

\setcounter{table}{5}
\begin{table}[!ht]
\caption{Good $3$-element antichains in type $H_4$, continued}
\begin{tabular}{|c|r@{\,+\,$c$\,}l|}
\hline
  $I_{\min}$ &
\multicolumn{2}{c|}{
  Solutions to $(v \mid \beta)=1$ for all $\beta \in I_{\min}$}
\\ \hline\hline

   $\{\alpha_{41},\alpha_{42},\alpha_{43}\}$
 & $\frac{1}{4}(4\tau-6,\,2\tau-3,\,5-3\tau,\,5\tau-8)$&$(-2\tau,\,\tau,\,1,\,\tau-1)$
\\ \hline

   $\{\alpha_{42},\alpha_{43},\alpha_{44}\}$
 & $\frac{1}{4}(10-6\tau,\,5-3\tau,\,2-\tau,\,2\tau-3)$&$(2\tau,\,-\tau-2,\,1,\,\tau-1)$
\\ \hline

   $\{\alpha_{43},\alpha_{44},\alpha_{45}\}$
 & $\frac{1}{2}(5-3\tau,\,5\tau-8,\,5-3\tau,\,5-3\tau)$&$(1,\,\tau-3,\,1,\,3-2\tau)$
\\ \hline

   $\{\alpha_{44},\alpha_{45},\alpha_{46}\}$
 & $\frac{1}{4}(2\tau-3,\,10-6\tau,\,2\tau-3,\,10\tau-16)$&$(-\tau,\,2,\,-\tau,\,2\tau-2)$
\\ \hline
\end{tabular}
\end{table}

In Table~\ref{ta:comparisons} the comparison column for a given
$\Lambda$ is left blank when a prior comparison also applies to
yield that $R_\Lambda=\emptyset$.

\begin{table}[!ht]
\caption{Comparisons that reveal empty regions in type $H_4$}
\label{ta:comparisons}
\begin{tabular}{|c|c|r@{\;}l|}
\hline
  $I_{\min}$ & $I^c_{\max}$ &
\multicolumn{2}{c|}{
  Comparison of convex combinations}
\\ \hline\hline
 $\{\alpha_{16}\}$
 & $\{\alpha_{14}, \alpha_{25}\}$
 & $\alpha_{16}<$ &$\frac{1}{\tau+1}(\tau\alpha_{14}+\alpha_{25})$
\\ \hline

 $\{\alpha_{13}, \alpha_{16}\}$
 & $\{\alpha_{14}, \alpha_{19}\}$
 & $\alpha_{16}<$ &$\frac{1}{\tau+1}(\tau\alpha_{14}+\alpha_{19})$
\\ 

 $\{\alpha_{16}, \alpha_{17}\}$
 & $\{\alpha_{13}, \alpha_{14}, \alpha_{19}\}$
 &&
\\ 

 $\{\alpha_{16}, \alpha_{21}\}$
 & $\{\alpha_{14}, \alpha_{17}, \alpha_{19}\}$
 &&
\\ \hline

 $\{\alpha_{16}, \alpha_{25}\}$
 & $\{\alpha_{14}, \alpha_{21}\}$
 & $\alpha_{16}<$ &$\frac{1}{\tau+1}(\tau\alpha_{14}+\alpha_{21})$
\\ \hline

 $\{\alpha_{13}, \alpha_{20}\}$
 & $\{\alpha_{14}, \alpha_{16}, \alpha_{19}\}$
 & $\frac{1}{\tau+1}(\tau\alpha_{13}+\alpha_{20})<$ &$\frac{1}{\tau+1}(\tau\alpha_{14}+\alpha_{19})$
\\ \hline

 $\{\alpha_{14}, \alpha_{19}\}$
 & $\{\alpha_{24}\}$
 & $\frac{1}{\tau+1}(\tau\alpha_{14}+\alpha_{19})<$ &$\alpha_{24}$
\\ \hline

 $\{\alpha_{14}, \alpha_{21}\}$
 & $\{\alpha_{23}, \alpha_{24}\}$
 & $\frac{1}{\tau+1}(\tau\alpha_{14}+\alpha_{21})<$ &$\alpha_{24}$
\\ 

 $\{\alpha_{14}, \alpha_{21}, \alpha_{23}\}$
 & $\{\alpha_{19}, \alpha_{24}\}$
 &&
\\ \hline

 $\{\alpha_{14}, \alpha_{25}\}$
 & $\{\alpha_{28}\}$
 & $\frac{1}{\tau+1}(\tau\alpha_{14}+\alpha_{25})<$ &$\alpha_{28}$
\\ \hline

 $\{\alpha_{18}, \alpha_{25}\}$
 & $\{\alpha_{14}, \alpha_{28}\}$
 & $\frac{1}{\tau+1}(\tau\alpha_{18}+\alpha_{25})<$ &$\alpha_{28}$
\\ \hline

 $\{\alpha_{18}, \alpha_{21}\}$
 & $\{\alpha_{14}, \alpha_{23}, \alpha_{24}\}$
 & $\frac{1}{\tau+1}(\tau\alpha_{18}+\alpha_{21})<$ &$\frac{1}{\tau+1}(\alpha_{23}+\tau\alpha_{24})$
\\ \hline

 $\{\alpha_{21}, \alpha_{22}\}$
 & $\{\alpha_{18}, \alpha_{23}, \alpha_{24}\}$
 & $\frac{1}{\tau+1}(\alpha_{21}+\tau\alpha_{22})<$ &$\frac{1}{\tau+1}(\alpha_{23}+\tau\alpha_{24})$
\\ \hline

 $\{\alpha_{14}, \alpha_{25}, \alpha_{28}\}$
 & $\{\alpha_{21}, \alpha_{23}, \alpha_{24}\}$
 & $\frac{1}{\tau+1}(\alpha_{14}+\tau\alpha_{25})<$ &$\frac{1}{\tau+1}(\alpha_{21}+\tau\alpha_{24})$
\\ 

 $\{\alpha_{14}, \alpha_{23}, \alpha_{25}\}$
 & $\{\alpha_{21}, \alpha_{24}\}$
 &&
\\ \hline

 $\{\alpha_{18}, \alpha_{25}, \alpha_{28}\}$
 & $\{\alpha_{14}, \alpha_{21}, \alpha_{23}, \alpha_{24}\}$
 & $\frac{1}{\tau+1}(\alpha_{18}+\tau\alpha_{25})<$ &$\frac{1}{\tau+1}(\alpha_{23}+\tau\alpha_{24})$
\\ \hline
\end{tabular}
\end{table}

In Tables~\ref{ta:solutionsa} and~\ref{ta:solutionsb}, the values
given for $c$ and $d$ will cause the corresponding solution to lie
in $C$.  The solution will also satisfy $(v \mid \gamma)<1$ for
$\gamma \in I^c_{\max}$ in Table~\ref{ta:solutionsa} or $(v \mid
\beta)>1$ for $\beta \in I_{\min}$ in Table~\ref{ta:solutionsb},
as in the hypotheses of parts (i) and (ii) of
Theorem~\ref{thm:reduction-thru-antichain} respectively.  In some
cases, $R_\Lambda$ can be proven to be non-empty by applying the
conclusion of Theorem~\ref{thm:reduction-thru-antichain}(i) in
Table~\ref{ta:solutionsa} or
Theorem~\ref{thm:reduction-thru-antichain}(ii) in
Table~\ref{ta:solutionsb}.  In these cases, in place of a
solution we list the relationship between relevant increasing sets.
\begin{table}[!ht]
\caption{Solutions that reveal non-empty regions in type $H_4$
using Theorem~\ref{thm:reduction-thru-antichain}(i)}
\label{ta:solutionsa}
\begin{tabular}{|c|c|@{\quad}l|c|c|}
\hline
  $I_{\min}$ & $I^c_{\max}$ &
  Solutions to $(v \mid \beta)=1$ for all $\beta \in I_{\min}$
  & $c$ & $d$
\\ \hline\hline
 $\{\alpha_{13},\alpha_{19}\}$
 & $\{\alpha_{14},\alpha_{20}\}$
 &
$\frac{1}{2}(2\tau-3,\,2\tau-3,\,2-\tau,\,0)$
 & & \\
& & $\quad+c(1,\,-1,\,0,\,0)+d(0,\,0,\,0,\,1)$ & $0$ & $\frac{5\tau-8}{4}$
\\ \hline
 $\{\alpha_{14},\alpha_{23}\}$
 & $\{\alpha_{24},\alpha_{25}\}$
 & $(0,\, 2\tau-3,\, \frac{2-\tau}{2},\, \frac{2-\tau}{2})$
 & &
\\
& & $\quad+c(1,\,2-2\tau,\,\frac{\tau-1}{2},\,\frac{\tau-1}{2})+d(0,\,0,\,1,\,-1)$
& $5-3\tau$ & $\frac{\tau-2}{2}$ \\
 $\{\alpha_{18},\alpha_{23}\}$
 & $\{\alpha_{14},\alpha_{24},\alpha_{25}\}$
 & $I(\{\alpha_{18},\alpha_{23}\})=I(\{\alpha_{14},\alpha_{23}\})\setminus \{\alpha_{14}\}$
 &&
\\
 $\{\alpha_{14},\alpha_{28}\}$
 & $\{\alpha_{23},\alpha_{24},\alpha_{25}\}$
 & $I(\{\alpha_{14},\alpha_{28}\})=I(\{\alpha_{14},\alpha_{23}\})\setminus \{\alpha_{23}\}$
 &&
\\
 $\{\alpha_{18},\alpha_{28}\}$
 & $\{\alpha_{14},\alpha_{23},\alpha_{24},\alpha_{25}\}$
 & $I(\{\alpha_{18},\alpha_{28}\})=I(\{\alpha_{14},\alpha_{23}\})\setminus \{\alpha_{14},\alpha_{23}\}$
 &&
\\ \hline
 $\{\alpha_{14},\alpha_{24}\}$
 & $\{\alpha_{23},\alpha_{25}\}$
 & $(2-\tau,\, 0,\, 0,\, 2-\tau)$ & &
\\
& & $\quad+c(-\tau,\,1,\,0,\,0)
+d(-\tau,\,0,\,1,\,0)$ & $\frac{2-\tau}{4}$ & $\frac{2-\tau}{4}$ \\
 $\{\alpha_{18},\alpha_{24}\}$
 & $\{\alpha_{14},\alpha_{23},\alpha_{25}\}$
 & $I(\{\alpha_{18},\alpha_{24}\})=I(\{\alpha_{14},\alpha_{24}\})\setminus \{\alpha_{14}\}$
 &&
\\ \hline
\end{tabular}
\end{table}
\begin{table}[!ht]
\caption{Solutions that reveal non-empty regions in type $H_4$
using Theorem~\ref{thm:reduction-thru-antichain}(ii)}
\label{ta:solutionsb}
\begin{tabular}{|c|c|@{\quad}l|c|c|}
\hline
  $I_{\min}$ & $I^c_{\max}$ &
  Solutions to $(v \mid \gamma)=1$ for all $\gamma \in I^c_{\max}$
  & $c$ & $d$
\\ \hline\hline
 $\{\alpha_{16},\alpha_{19}\}$
 & $\{\alpha_{14},\alpha_{17}\}$
 & $(0,\,0,\,\tau-1,\,0)$
 & & \\
& & $\quad+c(1,\,0,\,-1,\,2-\tau)+d(0,\,1,\,-\tau,\,\tau-1)$ & $\frac{7\tau-7}{10}$ & $\frac{2-\tau}{5}$
\\ \hline
 $\{\alpha_{21},\alpha_{24}\}$
 & $\{\alpha_{22},\alpha_{23}\}$
 & $(0,\,2-\tau,\,0,\,0)$
 & & \\
& & $\quad+c(\tau-1,\,-1,\,1,\,0)+d(\tau-1,\,-1,\,0,\,1)$ &
$\frac{2-\tau}{2}$ & $\frac{2\tau-3}{4}$
\\ \hline
 $\{\alpha_{18},\alpha_{24},\alpha_{25}\}$
 & $\{\alpha_{14},\alpha_{21},\alpha_{23}\}$
 & $(1-\tau,\,1,\,0,\,0)+c(2\tau+1,\,-2\tau-2,\,\tau,\,1)$
 & $\frac{12-7\tau}{4}$ &
\\
 $\{\alpha_{14},\alpha_{24},\alpha_{25}\}$
 & $\{\alpha_{21},\alpha_{23}\}$
 & $I(\{\alpha_{14},\alpha_{24},\alpha_{25}\})=I(\{\alpha_{18},\alpha_{24},\alpha_{25}\})
\cup \{\alpha_{14}\}$
 &&
\\ \hline
\end{tabular}
\end{table}

The final antichain not listed in
Tables~\ref{ta:comparisons}-\ref{ta:solutionsb} requires a slight
modification of the proof of
Theorem~\ref{thm:reduction-thru-antichain}(ii) and is detailed in
the proof of Theorem~\ref{thm:H4}.

\end{appendix}


\bibliographystyle{amsalpha}

\begin{thebibliography}{[KR2]}

\bibitem{AIM}
 Prepared by D.~Armstrong,
 Braid groups, clusters, and free probability: An outline from
 the AIM workshop, January 2005.

\bibitem{A1}
 C.~Athanasiadis,
 Characteristic polynomials of subspace arrangements and finite fields,
 {Adv.\ Math.}\ {122} (1996), no.~2, 193--233.

\bibitem{A2}
 \bysame,
 Generalized Catalan numbers, Weyl groups and arrangements of hyperplanes,
 {Bull.\ London Math.\ Soc.}\ {36} (2004), no.~3, 294--302.

\bibitem{Ba1}
 D.~Barbasch,
 Spherical dual for $p$-adic groups,
 Geometry and representation theory of real and $p$-adic groups (C\'{o}rdoba, 1995),
 1--19, Progr.\ Math.\, 158, Birkhäuser, Boston, MA, 1998.

\bibitem{Ba2}
 \bysame,
 Unitary spherical spectrum for classical groups (to appear).

\bibitem{BM1}
 D.~Barbasch and A.~Moy,
 A unitarity criterion for p-adic groups,
 {Invent.\ Math.}\ {98} (1989), 19--37.

\bibitem{BM2}
 \bysame,
 Reduction to real infinitesimal character in affine Hecke algebras,
 {J.~Amer.\ Math.\ Soc.}\ {6} (1993), no.~3, 611--635.

\bibitem{BM3}
 \bysame,
 Unitary spherical spectrum for $p$-adic classical groups,
 {Acta.\ Appl.\ Math.}\ {44} (1996), no.~1-2, 3--37.

\bibitem{Be}
 D.~Bessis,
 The dual braid monoid,
 {Ann.\ Sci.\ \'Ecole Norm.\ Sup.}\ {36} (2003), no.~5, 647--683.

\bibitem{Br}
 T.~Brady,
 A partial order on the symmetric group and new $K(\pi,1)$'s for the
 braid groups,
 {Adv.\ Math.}\ {161} (2001), no.~1, 20--40.

\bibitem{BW1}
 T.~Brady and C.~Watt,
 $K(\pi,1)$'s for Artin groups of finite type,
 In Proceedings of the Conference on Geometric and Combinatorial Group Theory,
 Part I, Haifa, 2000,
 {Geom.\ Dedicata}\ {94}, (2002), 225--250.

\bibitem{BW2}
 \bysame,
 Lattices in finite real reflection groups, preprint 2005,
 \texttt{arXiv:math.CO/0501502}.

\bibitem{CP}
 P.~Cellini and P.~Papi,
 ad-nilpotent ideals of a Borel subalgebra II,
 {J.~Algebra}\ {258} (2002), 112--121.

\bibitem{Ch1}
 F.~Chapoton,
 Enumerative properties of generalized associahedra
 {S\'em.\ Lothar.\ Combin.}\ {51} (2004), Art.\ B51b, 16 pp.\ (electronic).

\bibitem{Ch2}
 \bysame,
 Sur le nombre de r\'eflexions pleines dans les groupes de Coxeter finis,
 preprint, 2004, {Bull.\ Belg.\ Math.\ Soc.}\ Simon Stevin (to appear),

\bibitem{C}
 D.~Ciubotaru,
 The unitary $\mathbb{I}$-spherical dual for split $p$-adic groups of type $F_4$,
 {Represent.\ Theory}\ {9} (2005), 94--137.

\bibitem{FR1}
 S.~Fomin and N.~Reading,
 {Root systems and generalized associahedra},
 Lecture notes for the IAS/Park City Graduate
 Summer School in Geometric Combinatorics, 2004.

\bibitem{FR2}
 S.~Fomin and N.~Reading,
 Generalized cluster complexes and
 Coxeter combinatorics,
 {Int.\ Math.\ Res.\ Not.}\ {44}
 (2005), 2709--2757.

\bibitem{FZ}
 S.~Fomin and A.~Zelevinsky,
 Y-systems and generalized associahedra,
 {Math.~Ann.}\ {158} (2003), 977--1018.


\bibitem{Humphreys:Reflection Groups}
 J.~E.~Humphreys,
 {Reflection groups and Coxeter groups},
 Cambridge Studies in Advanced Mathematics, 29.
 Cambridge University Press, Cambridge, 1990.

\bibitem{K}
 G.~Kreweras,
 Sur les partitions non crois\'ees d'un cycle,
 {Discrete Math.}\ {1} (1972), 333--350.

\bibitem{KR}
 C.~Kriloff and A.~Ram,
 Representations of graded Hecke algebras,
 {Represent.\ Theory} {6} (2002), 31--69 (electronic).

\bibitem{M}
 J.~McCammond,
 Noncrossing partitions in surprising locations,
 to appear in
 {Amer.\ Math.\ Monthly.}

\bibitem{P}
 D.~Panyushev,
 ad-nilpotent ideals of a Borel subalgebra: generators and duality,
 {J.\ Algebra}\ {274} (2004), no.~2, 822--846.

\bibitem{Rea}
 N.~Reading,
 Clusters, Coxeter-sortable elements and noncrossing partitions,
 preprint, December 2005, \texttt{arXiv:math.CO/0507186}.

\bibitem{Rei}
 V.~Reiner,
 Noncrossing partitions for classical reflection groups,
 {Discrete Math.}\ {177} (1997), no.~1-3, 195--222.

\bibitem{Sh1}
 J.-Y.~Shi,
 Sign types corresponding to an affine Weyl group,
 {J.~London Math.\ Soc.}\,
 (2) {35} (1987), 56--74.

\bibitem{Sh2}
 \bysame,
 The number of $\oplus$-sign types,
 {Quart.\ J.\ Math., Oxford Series (2)}\
 {48} (1997), no.~189, 93--105.

\bibitem{S}
 E.~Sommers,
 $B$-stable ideals in the nilradical of a Borel subalgebra,
 {Canad.\ Math.\ Bull.}\ {48} (2005), no.~3, 460--472.

\bibitem{St}
 J.~Stembridge,
 {A Maple Package for Posets}, computer software package
 available at http://www.math.lsa.umich.edu/$\sim$jrs/maple.html.

\bibitem{V}
 D.~Vogan,
 Computing the unitary dual, notes at atlas.math.umd.edu/papers.

\end{thebibliography}

\end{document}